\documentclass[12pt]{article}
\usepackage{amsmath,latexsym,amssymb}
\usepackage[mathscr]{eucal}
\newtheorem{thm}{Theorem}[section]
\newtheorem{lem}[thm]{Lemma}

\newtheorem{df}[thm]{Definition}

\newcommand{\id}{\mathrm{id}}

\newcommand{\cC}{\mathscr{C}}
\newcommand{\cG}{\mathscr{G}}
\newcommand{\cH}{\mathscr{H}}
\newcommand{\cK}{\mathcal{K}}
\newcommand{\cM}{\mathscr{M}}

\newcommand{\cN}{\mathscr{N}}
\newcommand{\cP}{\mathscr{P}}
\newcommand{\cQ}{\mathscr{Q}}
\newcommand{\cR}{\mathscr{R}}
\newcommand{\cU}{\mathscr{U}}
\newcommand{\ual}{w_\alpha}
\newcommand{\ube}{w_\beta}

\newcommand{\tM}{\tilde{\cM}}
\newcommand{\tilu}{\tilde{u}}
\newcommand{\tal}{\tilde{\alpha}}
\newcommand{\tbe}{\tilde{\beta}}

\newcommand{\circs}{{\lower-.3ex\hbox{{$\scriptscriptstyle\circ$}}}\hspace{1pt}}

\newcommand{\secp}{\tilde{p}}
\newcommand{\sq}{\tilde{q}}
\newcommand{\secpq}{\widetilde{pq}}

\newcommand{\sqr}{\widetilde{qr}}
\newcommand{\secr}{\widetilde{r}}

\DeclareMathOperator{\Ad}{Ad}
\DeclareMathOperator{\Tr}{Tr}
\DeclareMathOperator{\md}{mod}
\DeclareMathOperator{\Aut}{Aut}
\DeclareMathOperator{\Int}{Int}
\DeclareMathOperator{\Cntr}{Cnt_{\mathrm{r}}}

\DeclareMathOperator{\Ker}{Ker}

\begin{document}

\title{Classification of outer actions of discrete amenable groupoids on injective factors}
\author{Toshihiko Masuda \\
 Graduate School of Mathematics, Kyushu University, \\
 744, Motooka, Nishi-ku, 
Fukuoka, 819-0395, JAPAN, \\
e-mail: masuda@math.kyushu-u.ac.jp}
\date{}

\maketitle

\begin{abstract}
 We classify outer actions (or $\mathscr{G}$-kernels) of  discrete amenable groupoids on injective factors. 
Our method based on unified approach for classification of discrete amenable groups actions, and 
cohomology reduction theorem of discrete amenable equivalence relations. 
We do not use Katayama-Takesaki type resolution group approach. 
\end{abstract}

\section{Introduction}\label{sec:intro}

In the theory of operator algebras, research of automorphism groups and group actions is one of the central 
subject. In particular, classification of automorphisms and group actions has been developed by many hands 
since seminal works of A. Connes. Here we focus on classification of $G$-kernels and groupoid actions.

Classification of group actions and $G$-kernels was first taken place by A. Connes \cite{Co-peri}, 
\cite{Con-auto}. Since then,  
$G$-kernels on injective, semifinite factors has been classified by V. F. R. Jones \cite{J-act} for finite groups, 
and A. Ocneanu \cite{Ocn-act} for general discrete amenable groups. 
As the generalization of these works,
Katayama-Takesaki
classified $G$-kernels  of discrete amenable groups
on injective factors of type III.
(In what follows, we use the term ``outer actions'' instead of ``$G$-kernels'' according to Katayama-Takesaki.)
Their method based on clarification of cohomological 
aspect of outer actions, and they reduced classification to that of genuine group actions, 
which have been completely classified in \cite{KtST}.   

On the other hand, classification of discrete amenable measurable groupoids on injective factors has been also 
developed,  motivated by classification of compact abelian group actions.
In \cite{JT}, V. F. R. Jones and M. Takesaki classified actions of
compact abelian groups on semifinite injective factors. By the
Pontrjagin duality theorem, they first reduced classification to that of discrete
abelian group actions on semifinite von Neumann algebras. Then they reduced to
classification of groupoid actions on semifinite factors. 
This method  has been further extended by Sutherland-Takesaki
\cite{Su-Tak-RIMS}, and Kawahigashi-Takesaki \cite{Kw-Tak}.

In this paper, we further generalize these result and classify outer
action of an amenable discrete groupoid on injective factors.
We also realize actions with given invariant. 

Our main idea is similar to 
\cite{KtT-outerI}, \cite{KtT-outerII}, \cite{KtT-outerIII},
Namely, we split groupoid actions to isotropy part and  equivalence relation part.
We apply \cite{M-unif-Crelle} to classify isotropy part, and 
\cite{Su-coh-RIMS} to classify equivalence relation part. 
However there are two big difference between our argument and proceeding ones.

The first point is that we do not use the idea of resolution groups in
classification of outer actions. Indeed, Katayama and Takesaki's method
is the following. They clarified cohomological aspect of outer actions,
and reduced  all problems to those of genuine actions of discrete
amenable groups by means of resolution groups. 
One must note that the construction of resolution groups is not canonical. 
(We explain more detail in the beginning of \S \ref{sec:class}.)

The second point is that we do not use model action type argument in
our proof. In \cite{JT}, \cite{Su-Tak-RIMS}, \cite{Kw-Tak}, model actions
with special property were used to erase some obstruction arising in
classification. As explained above, we split classification to those of  
isotropy part and  equivalence relation part. Thus one must combine 
these two parts to obtain complete classification. In this procedure, some obstruction may arise to obtain 
cocycle conjugacy classification. 
In our approach, however, such obstruction does not arise, and we can avoid use of  model actions, 
since our classification of group actions in \cite{M-unif-Crelle} is slightly stronger than 
\cite{KtST} and  we  apply \cite{Su-coh-RIMS} more carefully.

This paper is organized as follows. In \S \ref{sec:preouter}, we give basic definition of outer actions of 
groupoids, introduce characteristic invariants of outer actions, and state main classification theorem. 
In \S \ref{sec:outergroup}, we show necessary facts on outer actions of a discrete amenable group.
In particular, we present the classification theorem of outer actions of discrete amenable groups on injective factors.
In \S \ref{sec:class}, we present the proof of main classification theorem. 
In \S \ref{sec:RWgroupoids}, we briefly review random walks of groupoids due to \cite{Kai-groupoid} \cite{ChoXin-Groupoid}, 
which is necessary to construct model actions.
In \S \ref{sec:model}, we construct model actions which realize given invariants. Main part is the construction of 
free actions on the injective factor of type II$_1$, where we use results in \S \ref{sec:RWgroupoids}.
In Appendix, we discuss relationship between our invariants and Katayama-Takesaki's invariants for outer actions
of discrete groups. 

\textbf{Acknowledgements.} The author is supported by 
JSPS KAKENHI Grant Number 16K05180. 

\section{Preliminaries on outer actions of  groupoids}\label{sec:preouter}

Our standard references for operator algebras are \cite{Tak-book}, and for
amenable groupoids is \cite{DeraGroupoidBook}.


Let $\cM$ be a von Neumann algebra, and $(\tM,\theta,\tau)$ be its continuous core
covariant system. Put $\cC=Z(\tM)$, and $\Aut_\theta(\cC):=\{\alpha\in
\Aut(\cC)\mid \alpha\circ\theta_t=\theta_t\circ \alpha, t\in \mathbb{R}\}$.
 Let $\tal\in \Aut(\tM)$ be the canonical
extension of $\alpha\in \Aut(\cM)$, the Connes-Takesaki module
$\md(\alpha)\in \Aut_\theta(\cC)$ the restriction of $\tal$ on $\cC$, and 
$\Cntr(\cM)=\{\alpha\in \Aut(\cM)\mid \tal \in \Int(\tM)\}$.

For a groupoid $\cG$, we use the following notation.
Its source map, range map and 
unit space are denoted by $s(g)$, $r(g)$ and $\cG^{(0)}$, respectively.
For $S\subset \cG$, $x\in \cG^{(0)}$,  $S^x:=\{g\in S\mid r(g)=x\}$,
$S_x:=\{g\in S \mid s(g)=x\}$,
$S_y^x=S_y\cap S^x$.
We denote  $x\sim y$, $x,y\in \cG^{(0)}$,  
if there is $g\in \cG$ with
$s(g)=x$, $r(g)=y$, which is an equivalence relation. 
Let
\[
 \cG^{(n)}:=\{(t_1,t_2,\cdots t_n)\mid t_i\in \cG,
s(t_i)=r(t_{i+1}), 1\leq i\leq n-1\}.
\]

Let $(\cG,\lambda,\nu)$ be 
a discrete measured groupoid.
Namely, $\cG^x$ is at most countable, $\nu$ is a quasi-invariant measure
on $\cG^{(0)}$, 
and $\lambda=\{\lambda^x\}_{x\in \cG^{(0)}}$ is the Haar system, 
where $\lambda^x$  is  the  countable measure on $\cG^x$.  
In the following, we simply write $(\cG,\lambda,\nu)$ by $\cG$.

Let $\cH$ and $\cK$ be its 
isotropy part  and equivalence relation part, i.e., 
$\cH=\{g\in \cG\mid s(g)=r(g)\}$, and 
$\cK=\{(y,x)\in G^{(0)}\times G^{(0)}\mid y\sim x \}$.
Hereafter we assume that $\cG$ is amenable, and ergodic.
Amenability of $\cG$ implies that amenability of $\cH_x$, and AF
property of $\cK$. Since $\cK$ is generated by a single ergodic transformation \cite{CoFeWei},
we can embed $\cK$ in $\cG$, and express $\cG$ as a semidirect product
$\cG=\cH\rtimes \cK$ (see \cite[p.1091]{Su-Tak-RIMS}).

For a measurable family of groups $\{G(x)\}_{x\in \cG^{0}}$,
let  $C^n(\cG,\{G(x)\})$ be a set of maps on $\cG^{(n)}$ such that 
$f(g_1,\cdots g_n)\in G(r(g_1))$, and $f(g_1,\cdots, g_n)=e$ if any of $g_i$ is in $\cG^{(0)}$.

\begin{df}\label{df:outeraction}
$(1)$ An outer action of $\cG$ (or $\cG$-kernel) on a measurable family of factors
 $\{\cM(x)\}_{x\in \cG^{(0)}}$ is a triple $(\alpha,w_\alpha,c)$ 
consisting of 
a map 
\[
 \alpha: g \in \cG\rightarrow \alpha_g\in \mathrm{Iso}\bigl(\cM(s(g)),\cM(r(g))\bigr), 
\]
$\ual\in C^2(\cG, \{\cU(\cM(x))\})$ and  $c\in C^3(\cG, \mathbb{T})$ 
satisfying
\begin{align*}
& \alpha_g\circ \alpha_h=\Ad \ual(g,h)\circ \alpha_{gh}, \\ 
& \alpha_g(w_\alpha(h,k))\ual(g,hk)=c(g,h,k)\ual(g,h)w_\alpha(gh,k).
\end{align*}
 
 Here $\mathrm{Iso}(\cM_1,\cM_2)$ is a set of all isomorphisms from $\cM_1$ onto $\cM_2$. 
Note that $c(g,h,k)$ satisfies the 3-cocycle identity
\[
 c(g,h,k) \overline{c(g,h,kl)} c(g,hk,l) \overline{c(gh,k,l)} c(h,k,l)=1,\,\,\,(g,h,k,l)\in \cG^{(4)}.
\]
$(2)$ Two outer actions $(\alpha,w_\alpha,c)$ and $(\beta,w_\beta,c)$ of
 $\cG$ with same 3-cocycle $c$ are said to be cocycle conjugate if there exist a
 measurable family 
 $\{\theta_x\}_{x\in \cG^{(0)}}$, $\theta_x\in \Aut(\cM(x))$, and 
 $u\in C^1(\cG, \{\cU(\cM(x))\})$ such that 
\[
 \Ad u(g)\circ \alpha_g=\theta_{r(g)}\circ \beta_g\circ \theta_{s(g)}^{-1},\,\,\,
 u(g)\alpha_g(u(h))w_{\alpha}(g,h)u(gh)^*=\theta_{r(g)}(w_\beta(g,h)).
\]
If we can take $\theta_x\in \overline{\Int}(\cM(x))$, we say that they are strongly
 cocycle conjugate.
\end{df}

\medskip

We introduce  invariants of $\alpha$ in a similar way as in the usual
group action case. 
Let $\{\tM(x),\theta^x\}$ and $\{\cC(x),\theta^x\}$ be the continuous core of $\cM(x)$, and the
flow of weights of $\cM(x)$, respectively.
The first invariant is the Connes-Takesaki module $\md(\alpha_g)$. 

We next introduce the  characteristic invariant of $(\alpha,w_\alpha,c)$.
Let $(\alpha^{x},\ual^x,c )$ be the  restriction of
$(\alpha,\ual,c)$ on $\cH_x$.
Of course $(\alpha^{x},\ual^x,c )$ is an outer action of $\cH_x$ on $\cM(x)$.
Let $\cN_{\alpha,x}\subset \cH_x$ be 
a normal subgroup defined by 
\[
 \cN_{\alpha,x}:=\{g\in \cH_{x}\mid \alpha_g^x\in \Cntr(\cM(x))\}.
\] 

Let
$V_\alpha(n,g):=\ual(g,g^{-1}n g)\ual(n,g)^*$ for 
$n, g\in \cG$. 
Then we have 
\[
 \alpha_g\circ \alpha_{g^{-1}ng}\circ
 \alpha_{g}^{-1}=\Ad (V_\alpha(n,g))\circ \alpha_{n}.
\]
Fix $\tilu_\alpha^x\in C^1\left(\cN_{\alpha,x},\cU(\tM(x))\right)$ such that $\tal_n^x=\Ad
\tilu_\alpha^x(n)$. Define $\lambda(n,g), \mu(m,n), d(n,t)\in \cU(\cC(x))$ by
\[
 \alpha_g(\tilu^{s(g)}(g^{-1}ng))=V_\alpha(n,g)\lambda(n,g)\tilu^{r(g)}(n),
\]
\[
\tilu_\alpha^x(m)\tilu_\alpha^x(n)=\mu(m,n)\ual^x(m,n)\tilu_\alpha^x(mn),
\]
\[
 \theta_t^x(\tilu_\alpha^x(n))=d(n,t)\tilu_\alpha^x(n).
\]

These unitaries enjoy the
following relations:
\begin{align*}
\mathrm{ (CC1)}& \,\,  \lambda(n,g)^*\theta_t(\lambda(n,g))=
d(n,t)^*\tal_g(d(g^{-1}ng,t)), \\
\mathrm{(CC2)}& \,\, d(m,t)d(n,t)d(mn,t)^* =\mu(m,n)^*\theta_t(\mu(m,n)), \\
\mathrm{(CC3) }& \,\, \mu(l,m)\mu(lm,n)=c(l,m,n)
\mu(m,n)\mu(l,mn), \\
\mathrm{(CC4)}&\,\, \lambda(n,gh)=\tal_g(\lambda(g^{-1}ng,h))\lambda(n,g) \\
&\hspace{52pt}\times \overline{c(g,g^{-1}ng, h)}
c(n,g,h)c(g,h,h^{-1}g^{-1}ngh), \\
\mathrm{(CC5)}&\,\, \lambda(mn,g)\lambda(m,g)^*\lambda(n,g)^*=
\mu(m,n)\alpha_g(\mu(g^{-1}mg,g^{-1}ng)^*)\\ 
&\hspace{52pt}\times {c(m,g,g^{-1}ng)} 
\overline{c(g,g^{-1}mg, g^{-1}ng)c(m,n,g)}, \\
\mathrm{(CC6)}&\,\, \lambda(n,m)=\mu(m,m^{-1}nm)\mu(n,m)^*, \\
\mathrm{(CC7)}&\,\, d(n,t+s)=\theta_t(d(n,s))d(n,t),
\end{align*} 
for $g,h\in \cG$, $m,n\in \cN$, $t,s\in \mathbb{R}$.

This $(\lambda,\mu,d)$ is called the characteristic cocycle of
$\alpha$. The triple $(\lambda,\mu,d)$ depends on the choice of
$\tilu^x(n)$. If we replace $\tilu^x(n)$ by $z^x(n)\tilu^x(n)$ for
$z^x(n)\in \cU(\cC(x))$, then the characteristic cocycle changes  into
\[
\left(\alpha_g(z^{s(g)}(g^{-1}ng))\lambda(n,g) z^{r(g)}(n)^*,
z^x(m)z^x(n)\mu(m,n)z^x(mn)^*, 
\theta_t(z^x(n))d(n,t)z^x(n)^*\right).
\]
The above relation is an equivalence relation on the set of characteristic cocycle for $(\alpha,w_\alpha,c)$.
Thus the equivalence class, say $\chi(\alpha)=[\lambda,\mu,d]$,  is a true invariant of
$\alpha$. 
Now we define the invariant of $\alpha$ by  $\mathrm{Inv}(\alpha)=(\cN_\alpha, \md(\alpha),\chi(\alpha))$.

The main classification theorem is as follows.
\begin{thm}\label{thm:classgroupoid}
Let $\cG$  an ergodic,  amenable 
 discrete  groupoid, 
$\{\cM(x)\}_{x\in \cG^{(0)}}$ be a measurable family of injective factors, 
and
 $(\alpha,\ual,c)$ and $(\beta,\ube,c)$ be two free outer actions 
 of $\cG$ on $\{\cM(x)\}_{x\in \cG^{(0)}}$ with same 3-cocycle $c$. Then $\alpha$ and
 $\beta$ are strongly  cocycle conjugate if and only if $\mathrm{Inv}(\alpha)=\mathrm{Inv}(\beta)$.
\end{thm}
The proof of Theorem \ref{thm:classgroupoid} will be presented in \S \ref{sec:class}.

\section{Preparation on outer actions of a group}\label{sec:outergroup}
In this section, we mainly treat group actions.

Let $G$ be a discrete group, and $\cM$ be a factor.
We gift the topology on $C^1(G, \cU(\cM))$ by
pointwise convergence in the strong topology of $\cU(\cM)$.

In a similar way as in \cite{M-unif-Crelle}, we can classify outer
actions of a discrete amenable group on an injective factor.

\begin{thm}\label{thm:classgroupouter}
 Let $\cM$ be an injective factor, $G$ a discrete amenable group.
Let $(\alpha,\ual,c)$, and $(\beta,\ube,c)$ be two outer actions of
 $G$ on $\cM$ with same 3-cocycle $c$. Then $\alpha$ and $\beta$ are
 strongly cocycle 
 conjugate if and only if $(N_\alpha,
 \md(\alpha),\chi(\alpha))=(N_\beta,\md(\beta),\chi(\beta))$ holds.
Namely, there exist
$v\in  C^1(G,\cU(\cM))$,
and $\theta \in \overline{\Int}(\cM)$ such that 
\[
 \Ad v(g)\circ \alpha_g=\theta\circ \beta_g\circ \theta^{-1}, 
v(g)\alpha_g(v(h))\ual(g,h)v(gh)^*=\theta(\ube(g,h)).
\]
\end{thm}
Since the proof is parallel to that of \cite{M-unif-Crelle}, we omit the proof.


Let $(\alpha,\ual,c)$ and $(\beta,\ube,c)$ be as in Theorem
\ref{thm:classgroupouter}, $N:=N_\alpha=N_\beta$, $[\lambda,\mu,d]:=\chi(\alpha)=\chi(\beta)$,
and fix $\tilu_\alpha, \tilu_\beta\in C^1(N, \cU(\tM))$ so that 
\begin{align*}
 \tal_g(\tilu_\alpha(g^{-1}ng))=\lambda(n,g)V_\alpha(n,g)\tilu_\alpha(n),  \,\,\,
& \tbe_g(\tilu_\beta(g^{-1}ng))=\lambda(n,g)V_\beta(n,g)\tilu_\beta(n), \\
\tilu_\alpha(m)\tilu_\alpha(n)=\mu(m,n)\ual(m,n)\tilu_\alpha(mn), \,\,\,
&\tilu_\beta(m)\tilu_\beta(n)=\mu(m,n)\ube(m,n)\tilu_\beta(mn), \\
\theta_t(\tilu_\alpha(n))=d(n,t)\tilu_\beta(n),\,\,\,
& \theta_t(\tilu_\alpha(n))=d(n,t)\tilu_\beta(n).
\end{align*}
If we carefully examine the proof 
in \cite{M-unif-Crelle}, 
it turns out that we can choose $v\in C^1(N, \cU(\cM))$ so that
\[
v(n)\tilu_\alpha(n)=\theta(\tilu_\beta(n)).
\]

This choice is crucial in our argument.

In what follows, we assume that $\cM$, $\alpha$ are as in Theorem
\ref{thm:classgroupouter}, and fix the choice of $\tilu_\alpha$.

Let $\Gamma_\alpha$ be a set of $(\theta,v)\in \Ker(\md)\times
C^1(G ,\cU(\cM))$
satisfying
\begin{align*}
 (\Gamma 1)\,\,\,& 
\Ad v(g)\circ \alpha_g=\theta\circ \alpha_g \circ \theta^{-1}, \,\,\, g\in G,\\
(\Gamma 2)\,\,\,&v(g)\alpha_g(v(h))\ual(g,h)v(gh)^*=\theta(\ual(g,h)),\,\,\, g,h\in G, \\
(\Gamma 3)\,\,\, &v (n)\tilu_\alpha(n)=\theta(\tilu_\alpha(n)),\,\,\, n\in N.
\end{align*}
We follow \cite{JT}, \cite{Su-Tak-RIMS}, \cite{Kw-Tak} for the
definition of $\Gamma_\alpha$. However, there is one different point
comparing with their definition, that is, 
we further assume the third condition. 
(Strictly speaking, our definition coincides with one in \cite{Su-Tak-RIMS}, 
and differ from ones in  \cite{JT}, \cite{Kw-Tak}.)
It turns out that the condition $(\Gamma 3)$ avoids us to use model actions in the proof of Theorem \ref{thm:classgroupoid}.

Define a group structure of $\Gamma_\alpha$ by 
\[
 (\theta_1,v_1)(\theta_2,v_2)=(\theta_1\theta_2,\theta_1(v_2(\cdot))v_1(\cdot)),\,\
(\theta,v)^{-1}=(\theta^{-1},\theta^{-1}(v(\cdot)^*)).
\]
By the product topology, we regard $\Gamma_\alpha$ as a Polish group.

Let $\Gamma_\alpha^0:=\{(\Ad w,w\alpha_g(w^*))\mid w\in \cU(\cM)\}$. Then
$\Gamma_\alpha^0$ is a normal subgroup of $\Gamma_\alpha$.

\smallskip

\noindent
\textbf{Remark.} When $\alpha$ is a genuine action of $G$,
$\Gamma_\alpha$ is identified with
\[
 \{\theta\in \Aut(\cM\rtimes_\alpha G)\mid \hat{\alpha}\circ
 \theta=(\theta\otimes \id)\circ \hat{\alpha}, \tilde{\theta}=\id
 \mbox{ on } \tM'\cap (\tM\rtimes_{\tal} G)\}
\]
as topological groups, where $\hat{\alpha}$ is the dual coaction, 
and 
$\Gamma^0_\alpha$ is $\{\Ad w\mid w\in \cU(\cM)\}$ via this identification.

\begin{lem}
We have $ \overline{\Gamma_\alpha^0}=\Gamma_\alpha$.
\end{lem}
\textbf{Proof.}
The proof is similar to that of \cite[Lemma 2.5.6]{JT}.
Fix a free ultrafilter $\omega$ over $\mathbb{N}$.
Take $(\theta,v)\in \Gamma_\alpha$. By \cite[Theorem 1]{KwST}, $\theta
\in \overline{\Int}(\cM)$.
Let $\{w_n\}\subset \cU(\cM)$ such that $\theta=\lim\limits_{n\rightarrow
\infty}\Ad w_n$, and set $W=(w_n)\in \cM^\omega$. 
Since $\theta\circ\alpha_g\circ \theta^{-1}=\Ad v(g)\circ \alpha_g$, $V(g):=W^*v(g)\alpha_g(W)\in \cM_\omega$.
We also have
$v(g)\alpha_g(v(h))\ual(g,h)v(gh)^*=\theta(\ual(g,h))=W\ual(g,h)W^*$. Thus 
\begin{align*}
 V(gh)\alpha_g(V(h))V(gh)^*
&=W^*v(g)\alpha_g(W)\alpha_g(W^*v(h)\alpha_h(W))V(gh)^* \\ 
 &=W^*v(g)\alpha_g(v(h))\ual(g,h)\alpha_{gh}(W)\ual(g,h)^*V(gh)^* \\ 
 &=W^*v(g)\alpha_g(v(h))\ual(g,h)\alpha_{gh}(W)\alpha_{gh}(W^*)v(gh)^*W\ual(g,h)^* \\ 
 &=W^*v(g)\alpha_g(v(h))\ual(g,h)v(gh)^*W\ual(g,h)^* \\ 
&=1.
\end{align*}
Hence $V(g)$ is a 1-cocycle for $\alpha_g$. 

We next show $V(n)=1$ for $n\in N_\alpha$. By the definition of
$\Gamma_\alpha$ and the fact $\cM^\omega\subset\tM^\omega$ \cite{Mato-disKac-Memo}, 
we have
\[
 v(n)\tilu_\alpha(n)=\theta(\tilu_\alpha(n))=W\tilu_\alpha(n)W^*,
\] and 
\[
 V(n)=W^*v(n)\alpha_n(W)=W^*v(n)\tilu_\alpha(n)W\tilu_\alpha(n)^*=1.
\]
Thus we can regard $V(\cdot)$ as a $G/N_\alpha$-cocycle on $\cM_\omega$. 
Again by \cite[Theorem 1]{KwST}, $\alpha$ is centrally free as an action of $G/N_\alpha$.
By the 1-cohomology vanishing \cite{Ocn-act}, we can
take $Z\in \cU(\cM_\omega)$
with $Z\alpha_g(Z^*)=V(g)$, and choose a representing sequence $WZ=(w'_n)$,
$w_n'\in \cU(\cM)$. Then $(\Ad w_n', w_n'\alpha_g(w_n^{'*}))\rightarrow
(\theta,v)$ as $n\rightarrow \omega$. \hfill $\Box$


\section{Classification}\label{sec:class}

In this section, we present a proof of Theorem \ref{thm:classgroupoid} by proving several lemmas.

Before the proof, we briefly explain Katayama-Takesaki's method in \cite{KtT-outerI}.
Assume a discrete amenable group $G$ 
and its 3-cocycle $c\in Z^3(G,\mathbb{T})$ is given. At first, they construct a discrete amenable group 
$G(c)$ (called the resolution group for $c$)
and its normal subgroup $H$ such that $G(c)/H=G$ and $[\pi_*(c)]\in H^3(G(c),\mathbb{T})$ is trivial, 
where $\pi:G(c)\rightarrow G$ is a quotient map. Then they reduce the problem to 
construction of genuine actions  of ${G}(c)$. Let $\alpha^c$  be an action  of $G(c)$ on a given  
injective factor $\cR$ such that $H\subset \{g\in G(c)\mid \alpha_g \in \Int(\cR)\}$. 
Fix a section $\mathfrak{s}:G=G(c)/H\rightarrow G(c)$. Then 
$\alpha_p:=\alpha^c_{\mathfrak{s}(p)}$ is an outer action of $G$ on $\cR$.
The modified HJR-exact sequence \cite[Theorem 2.7] {KtT-outerI} describes 
relation between $\mathrm{Inv}(\alpha)$ and $\mathrm{Inv}(\alpha^c)$, and
it is shown that any outer action of $G$  arises as above. 
In this way, they reduced all problems of outer actions to those of genuine actions, and obtained desired results by 
applying \cite{KtST}.   

Here we remark that $G(c)$ is not the canonical object, and heavily 
depends on the choice of representative of a 3-cocycle $c$ as remarked in \cite[Remark 2.15]{KtT-outerI}. 
This is an unsatisfactory point of Katayama-Takesaki's method.
Therefore, we do not take this approach, 
and classify and construct outer actions directly without use of resolution groups. 

In this section, we mainly denote elements of $\cH_x$ by roman alphabet $g,h,\cdots$, and 
those of $\cK$ by Greek letter $\gamma, \delta,\cdots$.
We can assume $\ual(\gamma,\delta)=\ube(\gamma,\delta)=1$ for $\gamma,\delta\in
\cK$, and hence $\alpha|_{\cK}$ and $\beta|_{\cK}$ are  genuine actions, 
since $\cK$ is generated by a single  transformation.

By Theorem \ref{thm:classgroupouter}, $\alpha^x$ and $\beta^x$ are
strongly cocycle conjugate for a.e. $x\in \cG^{(0)}$ with conjugating
automorphism $\theta_x$.
By replacing $\theta_{r(g)}\circ\beta_g\circ\theta_{s(g)}^{-1}$ with
$\beta_g$ and $\theta_x(\tilu_\beta^x(n))$ with $\tilu_\beta^x(n)$
respectively, 
we can assume that 
\begin{align*}
&  \Ad v^{x}(g)\circ \alpha_g^{x}=\beta_g^{x}, \\ 
 & \md(\alpha_\gamma)=\md(\beta_\gamma), \\ 
& v^{x}(g)\alpha_g^{x}(v^{x}(h))\ual^{x}(g,h)v^{x}(gh)^*=\ube^{x}(g,h), \\
& \alpha_\gamma(\tilu_\alpha^{x}(\gamma^{-1}n\gamma))=V_\alpha(n,\gamma)\lambda(n,\gamma)\tilu_\alpha^y(n), \\
 &\beta_\gamma(\tilu_\beta^{x}(\gamma^{-1}n\gamma))=V_\beta(n,\gamma)\lambda(n,\gamma)\tilu_\beta^y(n), \\
& v^x(m)\tilu_\alpha^x(m)=\tilu_\beta^x(m)
\end{align*}
 for $x\in \cG^{(0)}$, $g,h\in \cH_x$, $m\in \cN_x$, 
$n\in \cN_y$, $\gamma=(y,x)\in \cK$, and some $v^{x}\in C^1(\cH_x,\cU(\cM(x)))$.

Let 
$\Gamma(x):=\Gamma_{\alpha^{x}}$, $\Gamma^0(x):=\Gamma_{\alpha^{x}}^0$.
As in \cite{JT}, we can see that $\{\Gamma(x)\}$ is a Borel family of Polish groups.

For our purpose, it is convenient to extend the definition of $\Gamma(x)$ as follows.
For $\gamma=(y,x)\in \cK$, $\xi\in  C^2\left(\cH_y, \mathbb{T}\right)$, 
and $a\in C^1\left(\cN_y, \cU(\cC(y))\right)$, 
let $\Gamma(\gamma,\xi,a)$ be the set of $(\sigma,u)\in \mathrm{Iso}(\cM(x),\cM(y))
\times C^1\left(\cH_y,\cU(\cM(y))\right)$
such that 
\begin{align*}
(\Gamma 0)' &\,\,\,  \md(\sigma)=\md(\alpha_\gamma), \\
(\Gamma 1)' &\,\,\, \sigma\circ \alpha^x_{\gamma^{-1}g\gamma} \circ \sigma^{-1}=\Ad u(g)\circ
 \alpha_g^y, \,\,\, g\in \cH_y,\\
(\Gamma 2)'&\,\,\,u(g)\alpha_g^y(u(h))\ual^y(g,h)u(gh)^*=\xi(g,h)
\sigma(\ual^x(\gamma^{-1}g\gamma,\gamma^{-1}h\gamma)),\,\,\, g,h\in \cH_y \\
(\Gamma 3)'&\,\,\,a(n)u (n)\tilu^y_\alpha(n)=\sigma(\tilu^x_\alpha(\gamma^{-1}n\gamma)),\,\,\, n\in \cN_y.
\end{align*}
Of course, we have $\Gamma((x,x),1,1)=\Gamma(x)$ in this notation.

Define  product
$(\sigma_1,u_1)(\sigma_2,u_2)$, and  inverse $(\sigma,u)^{-1}$ by
\[
\left(\sigma_1\sigma_2, \sigma_1(u_2(\gamma^{-1}_1g\gamma_1))u_1(g)\right),\,\,\,
(\sigma,u)^{-1}=\left(\sigma^{-1},\sigma^{-1}(u(\gamma g\gamma^{-1}))^*\right)
\]
for $(\sigma,u)\in \Gamma(\gamma,\xi,a)$, 
$(\gamma_1,\gamma_2)\in \cK^{(2)}$, 
$(\sigma_1,u_1)\in \Gamma(\gamma_1,\xi_1,a_1)$, and  $(\sigma_2,u_2)\in \Gamma(\gamma_2,\xi_2,a_2)$.

For $\xi\in C^2\left(\cH_x, \mathbb{T}\right)$, 
$a\in C^1\left(\cN_x,\cU\left(\cC(x)\right)\right)$, and $\gamma=(y,x)\in \cK$, 
define $\gamma_*(\xi)\in C^2\left(\cH_y, \mathbb{T}\right)$, 
and $\gamma_*(a)\in C^1(\cN_y,\cU(\cC(y)))$ by 
\[
 \gamma_*(\xi)(g,h)=\xi(\gamma^{-1}g\gamma, \gamma^{-1}hg)\,\,\, \gamma_*(a)(n)=\md(\alpha_\gamma)(a(\gamma^{-1}n\gamma)).
\]

\begin{lem}\label{lem:Gammaproperty}
We have following. \\
$(1)$ Let  $(\gamma_1,\gamma_2)\in \cK^{(2)}$. 
Then 
$(\sigma_1,u_1)(\sigma_2,u_2)\in \Gamma(\gamma_1\gamma_2,\xi_1\gamma_{1*}(\xi_2),a_1\gamma_{1*}(a_2))$ holds
for $(\sigma_1,u_1)\in \Gamma(\gamma_1,\xi_1,a_1)$, and  $(\sigma_2,u_2)\in \Gamma(\gamma_2,\xi_2,a_2)$. \\
$(2)$ For $(\sigma,u)\in \Gamma(\gamma,\xi,a)$, 
$(\sigma,u)^{-1}\in 
\Gamma(\gamma^{-1},\gamma^{-1}_*(\xi)^*, \gamma^{-1}_*(a)^*)$.  \\
$(3)$ 
We have $(\sigma_1,u_1)(\theta,v)(\sigma_2,u_2)^{-1}\in \Gamma(y)$ 
for $\gamma=(y,x)$, 
$(\sigma_1,u_1),(\sigma_2,u_2)\in \Gamma(\gamma,\xi,a)$, and $(\theta,v)\in \Gamma(x)$,
\end{lem}
\textbf{Proof.}
It is routine to show (1) and (2). Then (3) follows from (1) and (2). 
(Note $\md(\sigma_1)=\md(\sigma_2)=\md(\alpha_\gamma)$.)
\hfill$\Box$

\begin{lem}\label{lem:Vcocycle}
Put
\begin{align*}
\zeta_\gamma(g,h)
&:=\overline{c(g,h,\gamma)}c(g,\gamma,\gamma^{-1}h\gamma)\overline{c(\gamma,
\gamma^{-1}g\gamma,\gamma^{-1}h\gamma)}, \,\,\, g,h \in \cH, \gamma\in \cK, \\
\eta_{\gamma_1,\gamma_2}(g)&:=\overline{c(\gamma_1,\gamma_1^{-1}g\gamma_1,\gamma_2)}c(\gamma_1,\gamma_2,\gamma^{-1}_2\gamma_1^{-1}g\gamma_1\gamma_2)\overline{c(g,\gamma_1,\gamma_2)}, \,\,\, g\in \cH, \gamma_1,\gamma_2\in \cK.
\end{align*}
 We have
\begin{align*}
(1) &\,\,\, V_\alpha(g,\gamma)\alpha_{g}^{y}(V_\alpha(h,\gamma))w^{y}_\alpha(g,h)V_\alpha(gh,\gamma)^*=
\zeta_\gamma(g,h)
\alpha_\gamma(\ual^{x}(\gamma^{-1}g\gamma,\gamma^{-1}h\gamma)). \\
(2) &\,\,\, \alpha_{\gamma_1}(V_\alpha(\gamma_1^{-1}g\gamma_1,\gamma_2))V_\alpha(g,\gamma_1)=
\eta_{\gamma_1,\gamma_2}(g)V_\alpha(g,\gamma_1\gamma_2).
\end{align*}

\end{lem}
\textbf{Proof.}
We only present the proof of (1);
\begin{align*}
 &V_\alpha(g,\gamma)\alpha_{g}^{y}(V_\alpha(h,\gamma))w^{y}_\alpha(g,h)V_\alpha(gh,\gamma)^* \\
&=\ual^{y}(\gamma,\gamma^{-1}g\gamma)
 \ual^{y}(g,\gamma)^*
\alpha_{g}^{y}(\ual^{y}(\gamma,\gamma^{-1}h\gamma)
\ual^{y}(h,\gamma)^*)w^{y}(g,h) \\
&\times w^{y}_\alpha(gh,\gamma)
\ual^{y}(\gamma,\gamma^{-1}gh\gamma)^*
 \\
&=\overline{c(g,h,\gamma)}
\ual^{y}(\gamma,\gamma^{-1}g\gamma)
 \ual^{y}(g,\gamma)^*
\alpha_{g}^{y}(\ual^{y}(\gamma,\gamma^{-1}h\gamma))\ual^{y}(g,h\gamma)
\ual^{y}(\gamma,\gamma^{-1}gh\gamma)^*
 \\
&=\overline{c(g,h,\gamma)}c(g,\gamma,\gamma^{-1}h\gamma)
\ual^{y}(\gamma,\gamma^{-1}g\gamma)
\ual^{y}(g\gamma,\gamma^{-1}h\gamma)
\ual^{y}(\gamma,\gamma^{-1}gh\gamma)^*
 \\
&=\overline{c(g,h,\gamma)}c(g,\gamma,\gamma^{-1}h\gamma)\overline{c(\gamma, \gamma^{-1}g\gamma,\gamma^{-1}h\gamma)}
\alpha_\gamma(\ual^{x}(\gamma^{-1}g\gamma,\gamma^{-1}h\gamma)).
\end{align*}
The second equation can be verified in a similar way.
\hfill$\Box$

\medskip

\begin{lem}\label{Borelfunctor}
For $\gamma=(y,x)\in\cK$, set
 $W(g,\gamma)=\beta_\gamma(v^{x}(\gamma^{-1}g\gamma)^{*})V_\beta(g,\gamma)v^{y}(g)$.
 Let $\pi_1(\gamma)=(\alpha_\gamma, V_\alpha(g,\gamma))$, and $\pi_2(\gamma)=(\beta_\gamma,W(g,\gamma))$. \\
$(1) $ $\pi_1(\gamma),\pi_2(\gamma)\in \Gamma(\gamma, \zeta_\gamma, \lambda(\cdot,\gamma))$. \\
$(2) $ $\pi_i(\gamma_1)\pi_i(\gamma_2)=(\id, \eta_{\gamma_1,\gamma_2})\pi_i(\gamma_1\gamma_2)$ for 
$(\gamma_1,\gamma_2)\in \cK^{(2)}$ and $i=1,2$.
\end{lem}
\textbf{Proof.}
(1) By the definition of $V_\alpha(g,\gamma)$, $\lambda(n,\gamma)$, $\zeta_\gamma$ 
and Lemma \ref{lem:Vcocycle}(1), the fact
$\pi_1(\gamma)\in \Gamma(\gamma,\zeta_\gamma,\lambda(\cdot, \gamma))$ follows. 

For $\pi_2(\gamma)$, 
by the definition of $W(g,\gamma)$, $\lambda(n,\gamma)$, $\zeta_\gamma$ 
and Lemma \ref{lem:Vcocycle}(1), 
it is clear that $\pi_2(\gamma)$ satisfies 
$(\Gamma 0)'$, $(\Gamma 1)'$, $(\Gamma 2)'$.
We verify $(\Gamma 3)'$. By the choice of $\tilu^x_\alpha(n)$,
$\tilu^x_\beta(n)$, and $v^x(n)$, we have
\begin{align*}
\beta_\gamma(u_\alpha^x(\gamma^{-1}n\gamma))&=
\beta_\gamma(v^x(\gamma^{-1}n\gamma)^*)\beta_\gamma(u_\beta^x(\gamma^{-1}n\gamma))\\
&=
\lambda(n,\gamma)
\beta_\gamma(v^x(\gamma^{-1}n\gamma)^*)V_\beta(n,\gamma)u_\beta^y(n)\\
&=
\lambda(n,\gamma)
W(n,\gamma)u_\alpha^y(n).
\end{align*}

The statement (2) follows from Lemma \ref{lem:Vcocycle}(2).
\hfill$\Box$

\begin{lem}\label{lem:Borelfunctor1}
For $\gamma=(y,x)\in \cK$, define 
$F_\gamma:\Gamma(x)\rightarrow \Gamma(y)$ by 
\[
 F_\gamma(\theta,v):=\pi_1(\gamma)(\theta,v)\pi_1(\gamma)^{-1}
\]
Then $(\Gamma(x),F_\gamma)$ gives a Borel functor in the sense of
{\upshape  \cite[Definition 4.1]{Su-coh-RIMS}} with $F_\gamma(\Gamma^0_x)=\Gamma^0_y$.
\end{lem}
\textbf{Proof.}
By Lemma \ref{lem:Gammaproperty} and Lemma  \ref{Borelfunctor}, $F_\gamma:\Gamma(x)\rightarrow \Gamma(y)$ is well-defined. 
It is clear that $F_\gamma(\Gamma^0_x)=\Gamma^0_y$.
Since $(\id,\eta_{\gamma_1,\gamma_2})$ commutes with all elements of $\Gamma(r(\gamma_1))$,
$F_{\gamma_1}F_{\gamma_2}=F_{\gamma_1\gamma_2}$ by Lemma \ref{Borelfunctor}(2). \hfill$\Box$

\begin{lem}\label{lem:groupoidhom}
Define $\rho(\gamma)\in  \Gamma(y)$, $\gamma=(y,x)$,  by 
$\rho(\gamma)=\pi_2(\gamma)\pi_1(\gamma)^{-1}$.
Then $\rho$ is a cocycle in the sense $\rho(\gamma_1\gamma_2)=\rho(\gamma_1)F_{\gamma_1}(\rho(\gamma_2))$.
\end{lem}
\textbf{Proof.}
By Lemma \ref{lem:Gammaproperty} and Lemma \ref{Borelfunctor}, $\rho(\gamma)\in \Gamma(y)$.
Since $(\id,\eta_{\gamma_1,\gamma_2})$ 
commutes with all element of $\Gamma(y)$,
\begin{align*}
\rho(\gamma_1)F_{\gamma_1}(\rho(\gamma_2))
&=\pi_2(\gamma_1)\pi_1(\gamma_1)^{-1}\pi_1(\gamma_1)\pi_2(\gamma_2)\pi_1(\gamma_2)^{-1}
\pi_1(\gamma_1)^{-1} \\
&=\pi_2(\gamma_1)\pi_2(\gamma_2)\pi_1(\gamma_2)^{-1}\pi_1(\gamma_1)^{-1} \\
&= (\id,\eta_{\gamma_1,\gamma_2}) 
\pi_2(\gamma_1\gamma_2)\pi_1(\gamma_1\gamma_2)^{-1} (\id,\eta_{\gamma_1,\gamma_2})^{-1} \\
&=\pi_2(\gamma_1\gamma_2)\pi_1(\gamma_1\gamma_2)^{-1} \\
&= \rho(\gamma_1\gamma_2).
\end{align*}
\hfill$\Box$

\medskip

\bigskip

\noindent
\textbf{Remark.}
We can not
show $\pi_2(\gamma)\in \Gamma(\gamma,\zeta_\gamma,\lambda(\cdot,n))$, and 
$\rho(\gamma)\in\Gamma(y)$ without the condition $v^x(n)\tilu_\alpha^x(n)=\tilu_\beta^x(n)$, 
which is missing in the proof of
\cite[Lemma 4.2]{Su-Tak-RIMS}. 

We will explain this point in detail. 
In our choice, $v(n)$ satisfies the following three equation simultaneously:
\[
 v^x(n)\tilu_\alpha^x(n)=\tilu_\beta^x(n),\]
\[
\alpha_\gamma(\tilu_\alpha^x(\gamma^{-1}n\gamma))=\lambda(n,\gamma)
\tilu_\alpha^y(n),\,\,
\beta_\gamma(\tilu_\beta^x(\gamma^{-1}n\gamma))=\lambda(n,\gamma)
\tilu_\beta^y(n).
\]

By the traditional model action splitting argument, 
we can not get the condition  $v^x(n)\tilu_\alpha^x(n)=\tilu_\beta^x(n)$, and 
only have $v^x(n)\tilu_\alpha^x(n)=c^x(n)\tilu_\beta^x(n)$ for some $c^x(n)\in \cC(x)$.

In the rest of this remark, we use notation in \cite{Su-Tak-RIMS}.
(Actions $\alpha_k,m_k$ in \cite{Su-Tak-RIMS} correspond to $\beta_\gamma$ and $\alpha_\gamma$ in our notation.)
In \cite[p.1109]{Su-Tak-RIMS}, 
they define a unitary $a(n)=v(n)u(n)$, claim $\lambda_\alpha=\lambda_m$, and conclude
$\alpha_k(a(k^{-1}nk))=\lambda(n,k)a(n)$.
(Unitaries
$a(n)$, $v(n)$, $u(n)$
 correspond to $\tilu^x_\beta(n)$, $v^x(n)$, $\tilu_\alpha^x(n)$ respectively,
 in our notation.)
However in their definition, $\lambda_\alpha=\lambda_m$ is not clear, and 
in fact we can only say that
$\lambda_\alpha(n,k)=c(k^{-1}nk)c(n)^*\lambda_m(n,k)$ for some $c(n)\in \mathbb{T}$, and 
the proof of \cite[Lemma 4.2(2)]{Su-Tak-RIMS} fails.
(Of course, this problem occurs only if $\cN_{\alpha}$ is non-trivial.)

We also remark that such problem does not arise in \cite{JT}, \cite{Kw-Tak}, because they treat abelian groups, 
use the fact that any character of a subgroup of an abelian group can be extended to whole group,
and apply the cohomology lemma 
\cite[Theorem 5.5]{Su-coh-RIMS}
to groups which are slightly different from our $\Gamma(x)$.
Anyway some model action type argument can not be avoided.

\medskip

We continue the proof of Theorem \ref{thm:classgroupoid}.
It is obvious that $\gamma=(y,x)\rightarrow (\id,1)\in \Gamma(y)$ is a cocycle.
By \cite[Theorem 5.5]{Su-coh-RIMS}, there exists $P(x)\in \Gamma(x)$ such that
\[
 (\id,1)\equiv 
P(y)\rho(\gamma)
F_\gamma(P(x))^{-1}\mod\Gamma^0(y),
\]
i.e., for each $\gamma=(y,x)\in \cK$,  there exists $U(\gamma)\in \cU(\cM(y))$ such that
\[
 (\Ad U(\gamma), U(\gamma)\alpha_g^y(U(\gamma)^*))=
P(y)\rho(\gamma)
F_\gamma(P(x))^{-1}.
\]
 
\begin{lem}\label{lem:unitaryconnect}
Let $P(x)=(\theta_x,z^{x}(g))$, and  
$U(g)=\theta_x(v^{x}(g))z^{x}(g)$ for $g\in \cH_x$ .
Then 
\begin{align*}
 &\Ad U(g)\circ \alpha_g^x=\theta_x\circ \beta_g^x\circ \theta_x^{-1}, \,\,\, g\in \cH_x, \\
 &\Ad U(\gamma)\circ \alpha_\gamma =\theta_y\circ \beta_\gamma\circ \theta_x^{-1},\,\,\, \gamma=(y,x)\in \cK, 
\end{align*}
and 
\begin{align*}
&\theta_y(\ube(g,\gamma)^*)U(g)\alpha_g^y(U(\gamma))\ual(g,\gamma) \\
&=
\theta_y(\ube(\gamma,\gamma^{-1}g\gamma)^*)U(\gamma)\alpha_\gamma(U(\gamma^{-1}g\gamma))
\ual(\gamma,\gamma^{-1}g\gamma),\,\,\, g\in \cH_y, \,\,\,\gamma=(y,x)\in \cK.
\end{align*}
%
hold.
\end{lem}
\textbf{Proof.}
Since $(\theta_x,z^x(g))\in \Gamma(x)$, 
\[
 \Ad U(g)\circ \alpha_g=\Ad \left(\theta_x(v^x(g))z^x(g)\right)\circ \alpha_g^x
=\Ad \theta_x(v^x(g)) \circ \theta_x\circ\alpha_g^x\circ \theta_x^{-1}=
\theta_x\circ\beta^x_g\circ \theta_x^{-1}.
\]

By \begin{align*}
 (\Ad U(\gamma), U(\gamma)\alpha_g^y(U(\gamma)^*))& =
P(y)\rho(\gamma)
F_\gamma(P(x))^{-1}\\
&=P(y)\pi_2(\gamma)\pi_1(\gamma)^{-1}\pi_1(\gamma)P(x)^{-1}\pi_1(\gamma)^{-1}.
\end{align*}
we have
\[
 \left(\Ad U(\gamma), U(\gamma)\alpha_g^y(U(\gamma)^*)\right)\pi_1(\gamma)P(x)
=P(y)\pi_2(\gamma).
\]
The left hand side is 
\begin{align*}
\lefteqn{ \left(\Ad U(\gamma), U(\gamma)\alpha_g^y(U(\gamma)^*)\right)\pi_1(\gamma)P(x)} \\
&=
 \left(\Ad U(\gamma), U(\gamma)\alpha_g^y(U(\gamma)^*)\right)(\alpha_\gamma, V_\alpha(g,\gamma))
(\theta_x,z^x(g)) \\
&=
\left(\Ad U(\gamma), U(\gamma)\alpha_g^y(U(\gamma)^*)\right)
\left(\alpha_\gamma\theta_x,\alpha_\gamma(z^x(\gamma^{-1}g\gamma))V_\alpha(g,\gamma)\right) \\
&=
\left(\Ad U(\gamma)\alpha_\gamma\theta_x,
U(\gamma)\alpha_g(z^x(\gamma^{-1}g\gamma))V_\alpha(g,\gamma)\alpha_g^y(U(\gamma))^*\right).
\end{align*}
The right hand side is 
\begin{align*}
 P(y)\pi_2(\gamma)&=(\theta_y,z^y(g))(\beta_\gamma,W(g,\gamma)) \\
&=\left(\theta_y\beta_\gamma,\theta_y(W(g,\gamma))z^y(g)\right).
\end{align*}
By comparing the first component, we have 
$\Ad U(\gamma)\circ \alpha_\gamma=\theta_y\circ \beta_\gamma\circ \theta_x^{-1}$. 

We next compare the second component;  
\begin{align*}
{U(\gamma)\alpha_g(z^x(\gamma^{-1}g\gamma))V_\alpha(g,\gamma)\alpha_g^y(U(\gamma))^*} 
&=\theta_y(W(g,\gamma))z^y(g) \\
&=\theta_y\beta_\gamma(v^x(\gamma^{-1}g\gamma))^*\theta_y(V_\beta(g,\gamma))\theta_y(v^y(g))z^y(g) \\
&=U(\gamma)\alpha_\gamma(v^x(\gamma^{-1}g\gamma))^*U(\gamma)^*\theta_y(V_\beta(g,\gamma))U(g).
\end{align*}
In the last equality, we used $\Ad U(\gamma)\circ \alpha_\gamma=\theta_y\circ \beta_\gamma\circ \theta_x^{-1}$.
Hence we obtain
\[
U(\gamma)
\alpha_\gamma(U(\gamma^{-1}g\gamma))
V_\alpha(g,\gamma)\alpha_g^{y}(U(\gamma))^*=\theta_y(V_\beta(g,\gamma))U(g).
\]
Hence we get the conclusion by $V_\alpha(g,\gamma)=\ual(\gamma,\gamma^{-1}g\gamma)\ual(g,\gamma)^*$.
\hfill$\Box$

\medskip

We extend the definition of $U(\cdot)$ for whole $\cG$ by 
\[
 U(g\gamma):=
\theta_y(\ube(g,\gamma)^*)U(g)\alpha_g^y(U(\gamma))\ual(g,\gamma), \,\,\,
 g\in \cH_y, \,\,\,\gamma=(y,x)\in \cK.
\]
By the above result, we also have
\[
 U(g\gamma)=
\theta_y(\ube(\gamma,\gamma^{-1}g\gamma)^*)U(\gamma)\alpha_\gamma(U(\gamma^{-1}g\gamma))
\ual(\gamma,\gamma^{-1}g\gamma).
\]
We can easily see $\Ad U(g)\circ
\alpha_{g}=\theta_{r(g)}\circ\beta_g\circ\theta_{s(g)}^{-1}$, $g\in \cG$.

\medskip

\noindent
\textbf{Remark.}
Here we emphasize that the condition 
 $v^x(n)\tilu_\alpha^x(n)=\tilu_\beta^x(n)$ 
is essential to deduce Lemma \ref{lem:unitaryconnect}. Without this condition, we can say only
weaker statement, i.e.,
\[
 U(g\gamma)=\varphi(g,\gamma)
\theta_y(\ube(\gamma,\gamma^{-1}g\gamma)^*)U(\gamma)\alpha_\gamma(U(\gamma^{-1}g\gamma))
\ual(\gamma,\gamma^{-1}g\gamma)
\]
for some $\varphi(g,\gamma)\in \mathbb{T}$, which is an obstruction for cocycle conjugacy of actions. 
In \cite{JT},
\cite{Su-Tak-RIMS}, \cite{Kw-Tak}, such obstruction appears, and   
model actions with special property are used 
to get rid of $\varphi(g,\gamma)$. 

\medskip

By verifying following lemma, we finish the proof of Theorem \ref{thm:classgroupoid}. 
\begin{lem}
 We have
\[
 U(g)\alpha_g(U(h))\ual(g,h)U(gh)^*=\theta_{r(g)}(\ube(g,h)),\,\,
 g,h\in \cG.
\]
\end{lem}
\textbf{Proof.} 
We denote $\theta_{r(g)}$ by $\theta$ for simplicity.
We first show lemma for $\ual(g\gamma,h)$ and
$\ual(g\gamma,\delta)$ for $g\in\cH_y$, $h\in \cH_{x}$,
$\gamma,\delta\in\cK$.
First note that 
\[
 \ual(g\gamma,\delta)=\overline{c(g,\gamma,\delta)}\ual(g,\gamma)^*\ual(g,\gamma\delta).
\]
Then we have 
\begin{align*}
 \lefteqn{U(g\gamma)\alpha_{g\gamma}(U(\delta))\ual(g\gamma,\delta)U(g\gamma\delta)^*} \\
&=\overline{c(g,\gamma,\delta)}
\theta_y(\ube(g,\gamma)^*)U(g)\alpha_g(U(\gamma))\ual(g,\gamma)\alpha_{g\gamma}(U(\delta))
\ual(g,\gamma)^*\ual(g,\gamma\delta) \\
&\times \ual(g,\gamma\delta)^*\alpha_g(U(\gamma\delta)^*)U(g)^*\theta(\ube(g,\gamma\delta)) \\
&=\overline{c(g,\gamma,\delta)}
\theta_y(\ube(g,\gamma)^*)U(g)\alpha_g(U(\gamma))\alpha_{g}\alpha_\gamma(U(\delta))
\alpha_g(U(\gamma\delta)^*)U(g)^*\theta(\ube(g,\gamma\delta)) \\
&=\overline{c(g,\gamma,\delta)}
\theta_y(\ube(g,\gamma)^*)U(g)\alpha_g(U(\gamma\delta))
\alpha_g(U(\gamma\delta)^*)U(g)^*\theta(\ube(g,\gamma\delta)) \\
&=\overline{c(g,\gamma,\delta)}
\theta_y(\ube(g,\gamma)^*\ube(g,\gamma\delta))\\ 
&=\theta(\ube(g\gamma,\delta)).
\end{align*}
In a similar way, we have
\begin{align*}
 \ual(g\gamma,h)& =\ual(\gamma\gamma^{-1}g\gamma, h) \\
& =\overline{c(\gamma, \gamma^{-1}g\gamma,h)}\ual(\gamma,\gamma^{-1}g\gamma)^*
\alpha_\gamma(\ual(\gamma^{-1}g\gamma, h ))\ual(\gamma,\gamma^{-1}g\gamma h).
\end{align*}
Thus
\begin{align*}
\lefteqn{ U(g\gamma)\alpha_{g\gamma}(U(h))\ual(g\gamma,h)U(g\gamma
 h)} \\
 &=
\theta(\ube(\gamma,\gamma^{-1}g\gamma)^*)U(\gamma)\alpha_\gamma(U(\gamma^{-1}g\gamma))
\ual(\gamma,\gamma^{-1}g\gamma)
\alpha_{g\gamma}(U(h)) \\ 
&\hspace{12pt} \times 
\overline{c(\gamma, \gamma^{-1}g\gamma,h)}\ual(\gamma,\gamma^{-1}g\gamma)^*
\alpha_\gamma(\ual(\gamma^{-1}g\gamma, h ))\ual(\gamma,\gamma^{-1}g\gamma h) \\
&\hspace{12pt}\times 
\ual(\gamma,\gamma^{-1}g\gamma h)^*
\alpha_\gamma(U(\gamma^{-1}g\gamma h))^*
U(\gamma)^*\theta(\ube(\gamma,\gamma^{-1}g\gamma h)) \\
&=\overline{c(\gamma, \gamma^{-1}g\gamma,h)}
\theta(\ube(\gamma,\gamma^{-1}g\gamma)^*)U(\gamma)\alpha_\gamma(U(\gamma^{-1}g\gamma))
\alpha_\gamma\alpha_{\gamma^{-1}g\gamma}(U(h)) \\ 
& \hspace{12pt}\times 
\alpha_\gamma\left(\ual(\gamma^{-1}g\gamma, h )\right)
\alpha_\gamma\left(U(\gamma^{-1}g\gamma h)\right)^*
U(\gamma)^*\theta(\ube(\gamma,\gamma^{-1}g\gamma h)) \\
&=\overline{c(\gamma, \gamma^{-1}g\gamma,h)}
\theta(\ube(\gamma,\gamma^{-1}g\gamma)^*)\theta\beta_\gamma\theta^{-1}\left(U(\gamma^{-1}g\gamma)
\alpha_{\gamma^{-1}g\gamma}(U(h))\ual(\gamma^{-1}g\gamma, h
 )u_{\gamma^{-1}g\gamma h}^*\right) \\
&\hspace{12pt}\times \theta(\ube(\gamma,\gamma^{-1}g\gamma h)) \\
&=\overline{c(\gamma, \gamma^{-1}g\gamma,h)}
\theta(\ube(\gamma,\gamma^{-1}g\gamma)^*)\theta\beta_\gamma\left(\ube(\gamma^{-1}g\gamma, h )\right)
\theta(\ube(\gamma,\gamma^{-1}g\gamma h)) \\
&=\theta(\ube(g\gamma,h)),
\end{align*}
and hence the statement holds for any $\ual(g,h)$, $g\in \cG$,
$h\in \cH$, or $h\in \cK$.

Finally, by using the equality
\[
 \ual(g,h\gamma)=c(g,h,\gamma)\alpha_g(\ual(h,\gamma))^*\ual(g,h)\ual(gh,\gamma),\,\,
g\in \cG, h\in \cH,\gamma\in \cK,
\]
we can show the statement in a similar way as follows;
\begin{align*}
\lefteqn{ U(g)\alpha_g(U(h\gamma))\ual(g,h\gamma)U(gh\gamma)^*} \\
&=U(g)\alpha_g\left(
\theta(\ube(h,\gamma)^*)U(h)\alpha_h(U(\gamma))\ual(h,\gamma)
\right) \\
&\hspace{12pt}\times c(g,h,\gamma)\alpha_g(\ual(h,\gamma))^*\ual(g,h)\ual(gh,\gamma)U(gh\gamma)^* \\
&=c(g,h,\gamma)
U(g)\alpha_g\theta(\ube(h,\gamma)^*)
\alpha_g(U(h))\alpha_g\alpha_h(U(\gamma))
\ual(g,h)\ual(gh,\gamma)U(gh\gamma)^* \\
&=c(g,h,\gamma)
\theta\beta_g(\ube(h,\gamma)^*)U(g)
\alpha_g(U(h))\ual(g,h)\alpha_{gh}(U(\gamma))\ual(gh,\gamma)U(gh\gamma)^* \\
&=c(g,h,\gamma)
\theta\beta_g(\ube(h,\gamma)^*)\theta(\ube(g,h)\ube(gh,\gamma)) \\
&=\theta \left(w_\beta(g,h\gamma)\right).
\end{align*}
\hfill$\Box$

\medskip 

Thus we have shown Theorem \ref{thm:classgroupoid}, 
which is a generalization of \cite{Su-Tak-RIMS} to outer action case.

%

\section{Random walks on groupoids}\label{sec:RWgroupoids}

In this section, we make preparation for construction of model actions. 
For this purpose, we collect necessary facts on random walks on groupoids.
Contents of this section is mainly based on
\cite{Kai-groupoid}, \cite{ChoXin-Groupoid}.

Let $f$ be a function on $\cG$. We often denote its restriction on $\cG^x$
by $f^x$. 
Let $f$ be a function on $\cG^x$, and $g\in \cG_x^y$. Then
$g\cdot f(a):=f(g^{-1}a)$ is a function on $\cG^y$.

\begin{df}\label{df:familyprobmeasure}
 Let $\mu$ be a measurable function on $\cG$. We say $\mu$ is a family of
 probability measure if each $\mu^x$ is a probability measure on
 $\cG^x$, and often denote it by $\mu=\{\mu^x\}_{x\in \cG^{(0)}}$.
\end{df}

Let $f\in \ell^1(\cG^x)$, and 
${\mu}=\{\mu^x\}_{x\in \cG^{(0)}}$ be a family of probability
measures on $\cG$ such that the support of $\mu^x$ is included in $\cG^x$.
Define $f*\mu\in \ell^1(\cG^x)$ by
\begin{align*}
 f*\mu(a)&=
\sum_{h\in \cG_{s(a)}}f(ah^{-1})\mu^{r(h)}(h) 
=\sum_{y\sim x}\sum_{h\in \cG^y_{s(a)}}{f(ah^{-1})}\mu^y(h) \\
&=\sum_{y\sim x}\sum_{h\in \cG^x_{y}}{f(h)}\mu^y(h^{-1}a)
=\sum_{y\sim x }\sum_{h\in \cG^x_{y}}{f(h)}h\cdot \mu^y(a).
\end{align*}

\begin{lem}\label{lem:L1norm}
 We have $\|f*\mu\|_1\leq \|f\|_1$ for $f\in \ell^1(\cG^x)$, where equality holds when $f$ is positive.
\end{lem}
\textbf{Proof.}
\begin{align*}
\|f*\mu\|_1&=\sum_{a\in \cG^x}|f*\mu(a) | 
=\sum_{z\sim x}\sum_{a\in \cG^x_z}|f*\mu(a) | \\
&=\sum_{z\sim x }\sum_{a\in \cG^x_z}
\left|\sum_{(y,x)\in \cK}\sum_{h\in
 \cG^y_{z}}{f(ah^{-1})}\mu^y(h)\right| \\
&\leq \sum_{z\sim x}\sum_{a\in \cG^x_z}
\sum_{(y,x)\in \cK}\sum_{h\in
 \cG^y_{z}}\left|{f(ah^{-1})}\right|\mu^y(h) \\
&= \sum_{z\sim x}\sum_{(y,x)\in \cK}\sum_{h\in \cG^y_{z}}
\sum_{a\in \cG^x_y}
\left|{f(a)}\right|\mu^y(h) \\
&= \sum_{(y,x)\in \cK}\sum_{h\in \cG^y}
\sum_{a\in \cG^x_y}
\left|{f(a)}\right|\mu^y(h) 
= \sum_{y\sim x}
\sum_{a\in \cG^x_y}
\left|{f(a)}\right| \\
&=\|f\|_1
\end{align*}
If $f$ is positive, the forth inequality becomes an equality. \hfill$\Box$

\medskip

Let $\mu=\{\mu^x\}$, and $\nu=\{\nu^x\}$ be families of probability
measures.  By Lemma \ref{lem:L1norm}, $\{\mu^x*\nu\}$ is also a family
of probability measures. Hence we can define a family of probability measure $\mu*\nu$ by $(\mu*\nu)^x :=\mu^x*\nu$.

Let $\mu$ be a family of probability measures on $\cG$.
For each $x\in \cG^{(0)}$, define a positive operator $P^x_\mu$ on
$\ell^\infty(\cG^{x})$ by 
\[
 P^x_\mu(f)(g)=\sum_{h\in \cG^x}\mu^{s(g)}(g^{-1}h)f(h)
=\sum_{h\in \cG^{s(g)}}f(gh)\mu^{s(g)}(h).
\]

\begin{lem}\label{lem:Markovcomposition}
$P_\mu^x$ is a unital positive operator, and 
 $P_\mu^x P_\nu^x=P_{\mu*\nu}^x$.
\end{lem}
\textbf{Proof.}
Positivity of $P_\mu^x$ is clear. We first show $P_\mu^x(1)=1$;

\[
 P_\mu^x(1)(g)=\sum_{h\in \cG^x}\mu^{s(g)}(g^{-1}h)
=\sum_{h\in \cG^{s(g)}}\mu^{s(g)}(h)=1.
\]

We next show $P_\mu^xP_\nu^x=P_{\mu*\nu}^x$.
Let $y=s(g)$. Then 
\begin{align*}
 P_\mu^x P_\nu^x(f)(g)&=\sum_{z\sim x}\sum_{h\in
 \cG^x_z}\mu^y(g^{-1}h)P_\nu(f)(h) \\
&=\sum_{z\sim x}\sum_{h\in \cG^x_z}\sum_{k\in \cG^x}\mu^y(g^{-1}h)\nu^z(h^{-1}k)f(k) \\
&=\sum_{k\in \cG^x}\left(\sum_{z\sim x}\sum_{h\in \cG^x_z}\mu^y(g^{-1}h)\nu^z(h^{-1}k)\right)f(k) \\
&=\sum_{k\in \cG^x}\left(\sum_{z\sim x}\sum_{l\in \cG^z_{s(k)}}\mu^y(g^{-1}kl^{-1})\nu^z(l)\right)f(k) \\
&=\sum_{k\in \cG^x}\mu^y*\nu(g^{-1}k)f(k) \\
&=P_{\mu*\nu}^x(f)(g).
\end{align*}
\hfill$\Box$

\medskip

\begin{lem}\label{lem:measureRW}
For $\theta\in \ell^1(\cG^x)$, define $\theta
P_\mu^x\in \ell^1(\cG^x)$ by 
\[
 \langle \theta P_\mu^x, f\rangle = \langle \theta, P_\mu^x f\rangle,
 \,\, f\in \ell^\infty(\cG^x).
\]
 Then $\theta P_\mu^x$ is given by $\theta*\mu$.
\end{lem}
\textbf{Proof.}
\begin{align*}
 \langle \theta, P_\mu^x f\rangle 
&= \sum_{y\sim x }\sum_{g\in \cG^x_y}\theta(g)P_\mu^x(f)(g) \\
&= \sum_{y\sim x}\sum_{g\in \cG^x_y}\sum_{h\in \cG^x}\theta(g)\mu^y(g^{-1}h)f(h) \\
&= \sum_{h\in \cG^x}\left(\sum_{y\sim x}\sum_{g\in \cG^x_y}\theta(g)\mu^y(g^{-1}h)\right)f(h) \\
&= \sum_{h\in \cG^x_z}\theta*\mu(h)f(h) \\
&= \langle \theta*\mu, f\rangle 
\end{align*}
\hfill$\Box$

\medskip
%
%

In \cite[Proposition 5.4]{ChoXin-Groupoid}, 
C- H. Chu and X. Li proved the following result 
by using the Reiter condition, which is
equivalent to the amenability of $\cG$ \cite[Proposition 3.2.14]{DeraGroupoidBook},  
and concluded  that an  amenable
groupoid is Liouville, i.e., the Poisson 
boundary of $P^x_\mu $ is trivial.

\begin{thm}\label{thm:groupoidReiter2}
 Let $\cG$ be an amenable, discrete  measured groupoid. Then there exists a family of
 probability measures
 $\mu=\{\mu^x\}$ such that 
$\lim\limits_n\|g\cdot \mu^{*n,s(g)}-\mu^{*n,r(g)}\|_{1}=0$ for
 a.e. $x\in \cG^{(0)}$ and 
 every $g\in \cG^x$, where $\mu^{*n}=\underbrace{\mu*\cdots *\mu}_{\mbox{$n$-times}}$.
\end{thm}
If we look their proof, we can see that $\mu$ can be chosen so that $\mathrm{supp}(\mu^x)=\cG^x$.
Thus we always assume $\mathrm{supp}(\mu^x)=\cG^x$, when we apply Theorem \ref{thm:groupoidReiter2}. 

For our purpose, we need stronger result than the Liouville property.
\begin{thm}\label{thm:groupoidsTailTrivial}
 Let $\cG$ be an amenable, discrete measured groupoid, and take $\mu=\{\mu^x\}$
 as in Theorem \ref{thm:groupoidReiter2}. Then the tail boundary of
 $P_\mu^x$ is trivial for a.e. $x\in \cG^{(0)}$.
\end{thm}
\textbf{Proof.} 
We denote $P_\mu^x$ by $P$ for simplicity.
Let $\{u_n\}\subset \ell^\infty(\cG^x)$ be a bounded harmonic sequence,
i.e., $P(u_{n+1})=u_n$, $n\in \mathbb{N}$, and $\sup_n\|u_n\|_\infty=M<\infty$.
Since the tail boundary of $P$ is identified with the space of all bounded harmonic 
sequences \cite{Kai-Markov}, 
we only have to show that $u_n$ is a constant function for any $n\in \mathbb{N}$.
Let $e\in \cG^{(0),x}$ a unit element.
On one hand, $P^{k}(u_{n+k})=u_n$, and hence we have 
\[
 \delta_gP^k(u_{n+k})-\delta_eP^k(u_{n+k})=
\langle \delta_g-\delta_e, u_n\rangle=
u_n(g)-u_n(e).
\]

On the other hand, 
we have
\begin{align*}
\left|\delta_gP^k(u_{n+k})-\delta_eP^k(u_{n+k})\right|&=
\left|\langle \delta_g*\mu^{*k}-\delta_e*\mu^{*k}, u_{n+k}\rangle\right| \\ 
&\leq 
\|\delta_g*\mu^{*k}-\delta_e*\mu^{*k}\|_1\|u_{n+k}\|_\infty \\ & \leq 
M\|g\cdot \mu^{*k,s(g)}-\mu^{*k,x}\|_{1} \\ & \rightarrow 0
\end{align*}
as $k\rightarrow \infty$ for all $g\in \cG^x$
by Lemma \ref{lem:Markovcomposition}, Lemma \ref{lem:measureRW} and the
choice of $\mu$. 
Hence $u_n$ is a constant function, and $u_n=P(u_{n+1})=u_{n+1}$. \hfill $\Box$

\medskip

\noindent
\textbf{Remark.}
(1) By \cite[Theorem 2.1]{Kai-Markov}, the tail boundary of $P=P_\mu^x$ is
trivial if and only if $\lim\limits_{n\rightarrow
\infty}\|\delta_eP^{n+d}-\nu P^n\|_1=0$ for any 
$d\in \mathbb{N}$, and  probability
measure $\nu$ on $\cG^x$ with $\nu \prec \delta_eP^d$.
We can prove Theorem 
\ref{thm:groupoidsTailTrivial} by verifying this condition. \\
(2) If we assume $\cG$ is an etale groupoid, then $\cG$ is topologically
amenable if and only if 
measurewise amenable, i.e., $(\cG,\nu)$ is amenable for every
quasi-invariant measure $\nu$ on $\cG^{(0)}$ by 
\cite[Theorem 3.3.7]{DeraGroupoidBook}. 
Then $\cG$ has a
Reiter condition by \cite[Corollary 3.3.8]{DeraGroupoidBook}, and we can
get rid of ``a.e.'' in Theorem 
\ref{thm:groupoidsTailTrivial}. 

%

\section{Construction of model actions}\label{sec:model}

In this section, we construct an outer action of $\cG$ with given invariant. 
As stated in \S \ref{sec:class}, we do not use resolution group method of Katayama-Takesaki \cite{KtT-outerI}.  

We first construct a free outer action on the injective factor of type II$_1$ in \S \ref{subsec:freeouter}.
This construction is inspired by ones  in  subfactor theory \cite{Po-amen}, \cite{HaYa-amenable}.

Then we realize  outer actions with given invariant in \S \ref{subsec:actionInv}   
by the same  method used in \cite{M-unif-Crelle}.

\subsection{Existence of free actions}\label{subsec:freeouter}
Let $\cG$ be a discrete measurable groupoid. (For instance, we do not assume the amenability of $\cG$.) 
We use the following notations;
\[
 \cG^{x,(n)}:=\{(t_1,t_2,\cdots t_n)\in \cG^{(n)}\mid t_1\in \cG^x\}, 
\]
\[
 \cG^{(n)}_y:=\{(t_1,t_2,\cdots t_n)\in \cG^{(n)}\mid 
t_n\in \cG_y\}
\]
and $\cG^{x,(n)}_y=\cG^{x,(n)}\cap \cG^{(n)}_y$.

Define $B_n(k)$ and $A_n(x)$, $k\in \cG$, $x\in \cG^{(0)}$, $n\geq 1$, by  
\[
 B_n(k):=B\left(\ell^2\left(\cG_{s(k)}^{r(k),(n)}\right)\right),\,\,\,
 A_n(x):=\bigoplus_{k\in \cG^{x}} B_n(k).
\]

Let us denote a matrix unit of $B_n(k)$ by
$e_{\boldsymbol{s},\boldsymbol{t}}$, $\boldsymbol{s},\boldsymbol{t}\in \cG^{r(k),(n)}_{s(k)}$.
For $g\in \cG$, and  $\boldsymbol{t}=(t_1,t_2,\cdots t_n)\in \cG^{s(g),(n)}$,
define $g\boldsymbol{t}=(gt_1,t_2,\cdots t_n)\in \cG^{r(g),(n)}$.

Fix a 3-cocycle $c$ of $\cG$.
Define $u_g^n(k)$, $w^n(g,h)(k)$, $I_n(k,l)$,  by 
\[
 u_g^n(k):=\sum_{\boldsymbol{t}\in \cG^{s(g),(n)}_{s(k)}}
c(g,t_1,t_1^{-1}g^{-1}k)e_{g\boldsymbol{t},\boldsymbol{t}} 
\in B\left(l^2(\cG^{s(g),(n)}_{s(k)}), l^2(\cG^{r(g),(n)}_{s(k)})\right), \,\, n\geq 1.
\]
\[
 w^n(g,h)(k):=\sum_{\boldsymbol{t}\in \cG^{s(h),(n)}_{s(k)}}
\overline{c(g,h,t_1)}e_{gh\boldsymbol{t},gh\boldsymbol{t}} \in B_n(k) ,\,\, n\geq 1,
\]
\[
 I_n(k,l):=\sum_{\boldsymbol{t}\in \cG^{r(k),(n-1)}_{s(k)}}
\overline{c(t_1,t_1^{-1}k, l)}e_{\boldsymbol{t},\boldsymbol{t}}\otimes e_{ll}\in B_n(kl),\,\, n\geq 2. 
\]
Of course, $u_g^n(k)$ and $w^n(g,h)(k)$ are unitaries, and $I_n(k,l)I_n(k,l)^*=I_n(k,l)^*I_n(k,l)^*=1\otimes e_{ll}$.

Define two injective homomorphisms $\alpha_g^n:A_n(s(g))\rightarrow A_n(r(g))$, and 
$\phi_n^x:A_n(x)\rightarrow A_{n+1}(x)$ by 
\[
 \alpha_g^n(a)(k):=\Ad u_g^n(k) \left(a(g^{-1}k)\right),
\]
and 
\[
 \phi_n^x(a)(k):=\sum_{l\in \cG_{s(k)}}I^{n+1}(kl^{-1},l)\left(
a(kl^{-1})\otimes 1\right)I^{n+1}(kl^{-1},l)^*
\]
respectively for $a=\bigoplus_{k\in \cG^x} a(k)$.


%

%

\begin{lem}\label{lem:unitaryrelation}
We have following relations. \\
$(1)$ $ u_g^n(k)u_h^n(g^{-1}k)= c(g,h, h^{-1}g^{-1}k)w^{n}(g,h)(k)u^n_{gh}(k)$. \\
$(2)$ $  c(g,g^{-1}k,l)
I_{n+1}(k,l)(u_g^n(k)\otimes 1)=u_g^{n+1}(kl)I_{n+1}(g^{-1}k,l)$.
 \end{lem}
\textbf{Proof.}
(1) It is shown as follows;
\begin{align*}
 u_g^n(k)u_h^n(g^{-1}k)&=
\left(\sum_{\boldsymbol{s}\in \cG^{s(g),(n)}_{s(k)}}
c(g,s_1,s_1^{-1}g^{-1}k)e_{g\boldsymbol{s},\boldsymbol{s}}\right)
\left(\sum_{\boldsymbol{t}\in \cG^{s(h),(n)}_{s(k)}}
c(h,t_1,t_1^{-1}h^{-1}g^{-1}k)e_{h\boldsymbol{t},\boldsymbol{t}}\right) \\
&=
\sum_{\boldsymbol{t}\in \cG^{s(h),(n)}_{s(k)}}
c(g,ht_1,t_1^{-1}h^{-1}g^{-1}k)c(h,t_1,t_1^{-1}h^{-1}g^{-1}k)
e_{gh\boldsymbol{t},\boldsymbol{t}}  \\
&=
\sum_{\boldsymbol{t}\in \cG^{s(h),(n)}_{s(k)}}
 \overline{c(g,h,t_1)}
 c(g,h, h^{-1}g^{-1}k)
c(gh,t_1,t_1^{-1}h^{-1}g^{-1}k)
e_{gh\boldsymbol{t},\boldsymbol{t}}  \\
&=
 c(g,h, h^{-1}g^{-1}k)w^{n}(g,h)(k)u^n_{gh}(k).
\end{align*}
Here we apply the 3-cocycle identity 
\begin{align*}
 &c(g,h,t_1)
\overline{ c(g,h, h^{-1}g^{-1}k)}
 c(g,ht_1,t_1^{-1}h^{-1}g^{-1}k)
\overline{ c(gh,t_1,t_1^{-1}h^{-1}g^{-1}k)}
 c(h,t_1,t_1^{-1}h^{-1}g^{-1}k) \\
&=1
\end{align*}
for $(g,h,t_1,t_1^{-1}h^{-1}g^{-1}k)$ at the third equality. 

\noindent
(2) The left hand side becomes as follows;
\begin{align*}
\lefteqn{ I_{n+1}(k,l)(u_g^n(k)\otimes 1)} \\
&=
\left(\sum_{\boldsymbol{s}\in \cG^{r(k),(n)}_{s(k)}}
\overline{c(s_1,s_1^{-1}k,l)}e_{\boldsymbol{s},\boldsymbol{s}}\otimes e_{ll}\right)
\left(\sum_{\boldsymbol{t}\in \cG^{s(g),(n)}_{s(k)}}
{c(g,t_1,t_1^{-1}gk)}e_{g\boldsymbol{t},\boldsymbol{t}}\otimes 1\right) \\
&=
\sum_{\boldsymbol{t}\in \cG^{s(g),(n)}_{s(k)}}
\overline{c(gt_1,t_1^{-1}g^{-1}k,l)}
{c(g,t_1,t_1^{-1}gk)}e_{g\boldsymbol{t},\boldsymbol{t}}\otimes e_{ll}.
\end{align*}

The right hand side is as follows;
\begin{align*}
\lefteqn{u_g^{n+1}(kl)I_{n+1}(g^{-1}k,l) }\\ &=
 \left(\sum_{\boldsymbol{s}\in \cG^{s(g),(n+1)}_{s(l)}}{c(g,s_1,s_1^{-1}g^{-1}kl)}e_{g\boldsymbol{s},\boldsymbol{s}}\right) 
\left(\sum_{\boldsymbol{t}\in \cG^{s(g),(n)}_{s(k)}}
\overline{c(t_1,t_1^{-1}g^{-1}k,l)}e_{\boldsymbol{t},\boldsymbol{t}}\otimes e_{ll}\right) \\
&=
\sum_{\boldsymbol{t}\in \cG^{s(g),(n)}_{s(k)}}
{c(g,t_1,t_1^{-1}g^{-1}kl)}\overline{c(t_1,t_1^{-1}g^{-1}k,l)}
e_{g\boldsymbol{t},\boldsymbol{t}}\otimes e_{ll}.
\end{align*}
By the 3-cocycle identity
\[
 c(g,t_1,t_1^{-1}g^{-1}k)
\overline{ c(g,t_1,t_1^{-1}g^{-1}kl)}
 c(g,g^{-1}k,l)
\overline{ c(gt_1,t_1^{-1}g^{-1}k,l)}
 c(t_1,t_1^{-1}g^{-1}k,l)=1
\]
for $(g,t_1,t_1^{-1}g^{-1}k,l)$, we have
\[
 c(g,g^{-1}k,l)
 I_{n+1}(k,l)(u_g(k)\otimes 1)=
u_g^{n+1}(kl) I_{n+1}(g^{-1}k).
\]
\hfill$\Box$

\begin{lem}\label{lem:modelinductivelimit}
We have  following. \\
$(1)$ $w^{n+1}(g,h)=\phi_{n}^{r(g)}(w^{n}(g,h))$, $(g,h)\in \cG^{(2)}$. \\
 $(2)$ $\phi_n^{r(g)}\circ \alpha_g^n=\alpha_g^{n+1}\circ \phi_n^{s(g)}$. \\
$(3)$ $\alpha_g^n\circ \alpha_h^n=\Ad w^n(g,h)\circ \alpha_{gh}$, $(g,h)\in \cG^{(2)}$. \\
$(4)$ $ \alpha_g^{(n)}(w^1(h,k))w^n(g,hk)=c(g,h,k)w^n(g,h)w^n(gh,k)$, $(g,h,k)\in \cG^{(3)}$. 
  \end{lem}
\textbf{Proof.} It is easy to see (1).

(2) On one hand, we have
\begin{align*}
\lefteqn{ \phi_n^{r(g)}\alpha_g^n(a)(k)} \\
&=\sum_{l\in \cG_{s(k)}}I^{n+1}(kl^{-1},l)\left(
\alpha_g^n(a)(kl^{-1})\otimes 1\right)I^{n+1}(kl^{-1},l)^* \\
&=\sum_{l\in \cG_{s(k)}}I^{n+1}(kl^{-1},l)\left(
u^n_g(kl^{-1})a(g^{-1}kl^{-1})u_g^n(kl^{-1})^*\otimes 1\right)I^{n+1}(kl^{-1},l)^* \\
&=\sum_{l\in \cG_{s(k)}}I^{n+1}(kl^{-1},l)\left(u^n_g(kl^{-1})\otimes 1\right)\left(a(g^{-1}kl^{-1})\otimes 1\right)
\left(u_g^n(kl^{-1})^*\otimes 1\right)I^{n+1}(kl^{-1},l)^*.
\end{align*}

On the other hand, 
\begin{align*}
\lefteqn{ \alpha_g^{n+1}\phi_n^{s(g)}(a)(k)=u_{g}^{n+1}(k)\phi_n^{s(g)}(a)(g^{-1}k)u_g^{n+1}(k)^* }\\
&=
\sum_{l\in \cG_{s(k)}}u_{g}^{n+1}(k)
I^{n+1}(g^{-1}kl^{-1},l)\left(a(g^{-1}kl^{-1})\otimes 1\right)I^{n+1}(g^{-1}kl^{-1},l)^*u_g^{n+1}(k)^*.
\end{align*}
By Lemma \ref{lem:unitaryrelation}, 
$ \phi_n^{r(g)}\alpha_g^n(a)=\alpha_g^{n+1}\phi_n^{s(g)}(a)$ 
holds.

(3)
By Lemma \ref{lem:unitaryrelation}, we have
\begin{align*}
 \alpha^n_g\alpha^n_h(a)(k)
&=\Ad u_g(k)\left(\alpha_h(a)(g^{-1}k)\right) \\ 
&=\Ad \left(u_g(k)u_h(g^{-1}k)\right)\left(a(h^{-1}g^{-1}k)\right) \\
&=\Ad \left(w^{(n)}(g,h)u^n_{gh}(k)\right)\left(a(h^{-1}g^{-1}k)\right) \\
&=\Ad \left(w^{(n)}(g,h)\right)\alpha_{gh}^n(a)(k)
\end{align*}
and the statement (3) holds.

(4)
At first, we consider the case $n=1$.
\begin{align*}
\lefteqn{ \alpha_g^{1}(w^1(h,k))(l)w^1(g,hk)(l)} \\
&=
\left(\sum_{t\in \cG^{s(k)}_{s(l)}}
\overline{c(h,k,t)}e_{ghkt,ghkt}\right)
\left(\sum_{s\in \cG^{s(k)}_{s(l)}}
\overline{c(g,hk,s)}e_{ghks,ghks}\right) \\
&=
\sum_{t\in \cG^{s(k)}_{s(l)}}
\overline{c(h,k,t)
c(g,hk,t)}e_{ghkt,ghkt} \\
&=
\sum_{t\in \cG^{s(k)}_{s(l)}}c(g,h,k)
\overline{c(g,h,kt)
c(gh,k,t)}e_{ghkt,ghkt}\,\,\, (\mbox{by the 3-cocycle identity}) \\
&=c(g,h,k)
\left(\sum_{s\in \cG^{s(k)}_{s(l)}}
\overline{c(g,h,s)}e_{s,s}\right)
\left(\sum_{t\in \cG^{s(k)}_{s(l)}}
\overline{c(gh,k,t)}e_{ghkt,ghkt}\right) \\
&=c(g,h,k)w^1(g,h)(l)w^1(gh,k)(l).
\end{align*}
By (1) and (2), we have $ \alpha_g^{n}(w^n(h,k))w^n(g,hk)=c(g,h,k)w^n(g,h)w^n(gh,k)$ for all $n\in \mathbb{N}$. 
\hfill$\Box$

\medskip

By Lemma \ref{lem:modelinductivelimit}, we get the model action in C$^*$-level.
\begin{thm}\label{thm:modelC*}
 Let $A(x):=\lim_{n}(A_n(x),\phi_n^x)$ be an inductive limit
 C$^*$-algebra. Then we can define an isomorphism 
$\alpha_g:A(x)\rightarrow A(y)$, $g\in \cG_x^y$, by 
$\alpha_g(a)=\alpha_g^n(a)$, $a\in A_n(x)$, and  $(\alpha_g,w(g,h))$ is an
 outer action of $\cG$ with 3-cocycle $c(\cdot,\cdot,\cdot)$.
\end{thm}

%

We will construct a suitable state $\psi_n^x$ on $A(x)$ for W$^*$-completion. 
Let $\mu=\{\mu^x\}$ be a family of probability measures on $\cG$ with $\mathrm{supp}(\mu^x)=\cG^x$.
Denote the non-normalized trace on $B_n(k)\left(=B\left(\ell^2\left(\cG_{s(k)}^{r(k),(n)}\right)\right)\right)$ by $\Tr_{n,k}$.

Let
\[
 \rho_n(k):=\sum_{\boldsymbol{t}\in \cG^{r(k),(n)}_{s(k)}}
\sum_{y\sim r(k)}
{\mu}^{r(k)}(kt_n^{-1}t_{n-1}^{-1}\cdots t_{1}^{-1}(r(k),y))\mu^{y}((y,r(k))t_1)\mu^{x_2}(t_2)\cdots \mu^{x_n}(t_n)
e_{\boldsymbol{t},\boldsymbol{t}},
\]
\[
 \psi_n^x(a):=\sum_{k\in \cG^x}\Tr_{n,k}(\rho_n(k)a(k)),\,\, a\in A_n(x), 
\]
\[
 E_n^x(a)(k):=\sum_{l\in \cG^{s(k)}}\mu^{r(l)}(l)\left(\id \otimes \Tr_{1,l}\right)\left(I^{n+1}(k,l)^*a(kl)I^{n+1}(k,l)\right), 
a\in A_{n+1}(x).
\]
(Note that we insert $(y,r(k))$ before $t_1$ in the definition of $\rho_n(k)$.)

\begin{lem}\label{lem:state}
 $(1)$ $\psi_n^x$ is a faithful state on $A_n(x)$. \\
$(2)$ 
$\psi_{n+1}^x(\phi_n^x(a)b)=\psi_n^x(aE_n^x(b))$, $a\in A_n(x)$, $b\in A_{n+1}(x)$. 
Thus $E_n^x$ is a $\psi_{n+1}^x$-preserving conditional expectation via identification 
$A_n(x)$ with $\phi_n^x(A_n(x))\subset A_{n+1}(x)$. \\
$(3)$ $\alpha_g^n\circ E_n^{s(g)}=E_n^{r(g)}\circ \alpha_g^{n+1}$.
\end{lem}
\textbf{Proof.}
(1)  Let $k\in \cG^x$.
Since we insert $(y,r(k))$ before $t_1$ in the definition of $\rho_n(k)$, we have 
\begin{align*}
\lefteqn{\Tr_{n,k}(\rho_n(k))} \\
&=\sum_{y \sim x}
\sum_{\boldsymbol{t}\in \cG^{r(k),(n)}_{s(k)}} 
{\mu}^{r(k)}(kt_n^{-1}t_{n-1}^{-1}\cdots
 t_{1}^{-1}(r(k),y))\mu^{y}((y,r(k))t_1)\mu^{r(t_2)}(t_2)\cdots \mu^{r(t_n)}(t_n) \\
&=\sum_{y\sim x}
\sum_{\boldsymbol{t}\in \cG^{y,(n)}_{s(k)}} 
{\mu}^{r(k)}(kt_n^{-1}t_{n-1}^{-1}\cdots
 t_{1}^{-1})\mu^{y}(t_1)\mu^{r(t_2)}(t_2)\cdots \mu^{r(t_n)}(t_n) \\
&=\mu^{*(n+1),x}(k).
\end{align*}
Hence $\psi_n^x(1)=\sum_{k\in \cG^x}\mu^{*(n+1),x}(k)=1$.
Since each $\rho_n(k)$ is a non-singular positive element, $\psi_n$ is a faithful normal state.

(2)
At first note the following relation.
\begin{align*}
\lefteqn{\rho_{n+1}(kl)(1\otimes e_{ll})} \\
=&
\sum_{\boldsymbol{t}\in \cG^{r(k),(n+1)}_{s(l)}}
\sum_{y\sim r(k)}
{\mu}^x(klt_{n+1}^{-1}t_{n}^{-1}\cdots t_{1}^{-1}(r(k),y))\mu^{y}((y,r(k))t_1)\mu^{x_2}(t_2)\cdots \mu^{x_{n+1}}(t_{n+1})
e_{\boldsymbol{t},\boldsymbol{t}} \\
 &\times (1\otimes e_{ll}) \\
=&
\sum_{\boldsymbol{t}\in \cG^{r(k),(n)}_{s(k)}}
\sum_{y\sim r(k)}
{\mu}^x(kt_{n}^{-1}\cdots t_{1}^{-1}k^{-1})\mu^{r(k)}(kt_1)\mu^{x_2}(t_2)\cdots \mu^{x_n}(t_n)\mu^{r(l)}(l)
e_{\boldsymbol{t},\boldsymbol{t}}\otimes e_{ll} \\
=&\mu^{r(l)}(l)\rho_{n}(k)\otimes e_{ll}.
\end{align*}

On one hand, we have 
\begin{align*}
\lefteqn{ \psi_{n+1}^x(\phi_n^x(a)b)} \\
&=
\sum_{k\in \cG^{x}}\Tr_{n+1,k}\left(\rho_{n+1}(k)\phi_n^x(a)(k)b(k)\right) \\
&=
\sum_{k\in \cG^x}\sum_{l\in \cG_{s(k)}}\Tr_{n+1,l}\left(\rho_{n+1}(k)I_{n+1}(kl^{-1},l)\left(a(kl^{-1})\otimes e_{ll}\right)I_{n+1}(kl^{-1},l)^*b(k)\right) \\
&=
\sum_{k\in \cG^x}\sum_{l\in \cG^{s(k)}}\Tr_{n+1,kl}\left(
\rho_{n+1}(kl)(1\otimes e_{ll})I_{n+1}(k,l)\left(a(k)\otimes 1\right)
I_{n+1}(k,l)^*b(kl)\right) \\ 
&=
\sum_{k\in \cG^x}\sum_{l\in \cG^{s(k)}}\mu^{r(l)}(l)\Tr_{n+1,kl}\left(
(\rho_{n}(kl)\otimes e_{ll})I_{n+1}(k,l)\left(a(k)\otimes 1\right)
I_{n+1}(k,l)^*b(kl)\right) \\
\end{align*}

On the other hand, 
\begin{align*}
 \psi_n^x(aE_n^x(b))&=\sum_{k\in \cG^x}\Tr_{n,k}\left(\rho_{n}(k)a(k)E_n^x(b)(k)\right) \\
&=\sum_{k\in\cG^x}\sum_{l\in \cG^{s(k)}}\mu^{r(k)}(l)\Tr_{n,k}\left(\rho_{n}(k)a(k)
(\id\otimes \Tr_{1,l})\left(I_{n+1}(k,l)^*b(kl)I_{n+1}(k,l)\right)
\right) \\
&=\sum_{k\in \cG^x}\sum_{l\in \cG^{s(k)}}
\mu^{r(l)}(l)\Tr_{n+1,kl}\left(\left(\rho_{n}(k)a(k)\otimes 1\right)
I_{n+1}(k,l)^*b(kl)I_{n+1}(k,l)\right) \\
&=\sum_{k\in \cG^x}\sum_{l\in \cG^{s(k)}}\mu^{r(l)}(l)\Tr_{n+1,kl}
\left(I_{n+1}(k,l)\left(\rho_{n}(k)a(k)\otimes 1\right)
I_{n+1}(k,l)^*b(kl)\right) \\
&=\sum_{k\in \cG^x}\sum_{l\in \cG^{s(k)}}\mu^{r(l)}(l)\Tr_{n+1,kl}\left((\rho_n(k)\otimes e_{ll})I_{n+1}(k,l)\left(a(k)\otimes 1\right)
I_{n+1}(k,l)^*b(kl)\right) .
\end{align*}
Here note $I_{n+1}(k,l)$ and $\rho_n(k)\otimes e_{ll}$ commute.
Hence we have 
$\psi_{n+1}^x(\phi_n^x(a)b)=\psi_n^x(aE_n^x(b))$, $a\in A_n(x)$, $b\in A_{n+1}(x)$, and 
$E_n^x$ is a $\psi_{n+1}^x$-preserving conditional expectation. 

(3) It is shown as follows;
\begin{align*}
\lefteqn{ E_n^{r(g)}\alpha_g^{n+1}(a)(k)} \\
&=\sum_{l\in \cG^{s(k)}}\mu^{r(l)}(l)
\left(\id\otimes\Tr_{1,l}\right)(I_{n+1}(k,l)^*\alpha_g^{n+1}(a)(kl)I_{n+1}(k,l)) \\
&=\sum_{l\in \cG^{s(k)}}\mu^{r(l)}(l)\left(\id\otimes 
\Tr_{1,l}\right)(I_{n+1}(k,l)^*u_g^{n+1}(kl)a(g^{-1}kl)u_g^{n+1}(kl)^*I_{n+1}(k,l)) \\
&=\sum_{l\in \cG^{s(k)}}\mu^{r(l)}(l)\left(\id\otimes 
\Tr_{1,l}\right)\left((u_g^n(k)\otimes 1)I_{n+1}(g^{-1}kl)a(g^{-1}kl)I_{n+1}(g^{-1}kl)(u_g^n(k)^*\otimes 1)\right) \\
&=\left(\id\otimes\Tr_{1,l}\right)
\left((u_g^n(k)\otimes 1)\phi_n^{s(g)}(a)(g^{-1}k)(u_g^n(k)^*\otimes 1)\right) \\
&=\alpha_g^nE_n^{s(g)}(a)(k).
\end{align*}
In the third equality, we used Lemma \ref{lem:unitaryrelation}. \hfill $\Box$

\medskip

By Lemma \ref{lem:state}, we can define a state $\psi^x$ on
$A(x)$ by $\psi^x(a)=\psi^x_n(a)$, $a\in A_n(x)$. Let $(\pi_x, H_x,\xi_x)$ be
the GNS triple of $A(x)$ by $\psi_x$, and 
$\cM(x):=\pi_x(A(x))''$.

Let $\sigma^{\psi^x}$ be the modular automorphism group. By Lemma \ref{lem:state}, 
the restriction of $\sigma^{\psi^x}$  on
$A_n(x)$ is $\sigma^{\psi^x_n}$, and there exists the $\psi^x$-preserving conditional expectation 
$F_n^x:\cM(x)\rightarrow A_n(x)$. 

\begin{thm}\label{thm:KMS-factor}
The center $Z(\cM(x))$ is identified with the tail boundary of $P^x_\mu$. 
Hence if $\mu$ satisfies the condition in Theorem \ref{thm:groupoidReiter2}, then $\cM(x)$ is a factor for a.e. $x$. 
\end{thm}
\textbf{Proof.} 
We simply denote $F_n^x$ and $E_n^x$ by $F_n$ and  $E_n$, respectively.
Take $z\in Z(\cM(x))$, and set $w_n=F_n(z)$. Then $w_n\in Z(A_n(x))=\ell^\infty(\cG^x)$, 
and $\|w_n\|_\infty \leq \|z\|$.  
Since $F_n=E_nF_{n+1}$, $E_{n}(w_{n+1})=w_n$, and 
\begin{align*}
 E_n(w_{n+1})(k)&=\sum_{l\in \cG^{s(k)}}\mu^l(l)(\id\otimes \Tr_{1,l})
(I(k,l)(w_{n+1}(kl)\otimes e_{ll})I_{n+1}(k,l)^*) \\
&=\sum_{l\in \cG^{s(k)}}w_{n+1}(kl)\mu^{r(l)}(l)(\id\otimes \Tr_{1,l})(1\otimes e_{ll})) \\
&=\sum_{l\in \cG^{s(k)}}w_{n+1}(kl)\mu^{r(l)}(l) \\
&=P^x_\mu(w_{n+1})(k)
\end{align*}
holds. Hence $\{w_n\}_n$ is a bounded harmonic sequence for $P_\mu^x$. Note that 
$\lim_nw_n=z$ in the $\sigma$-weak topology by the martingale convergence theorem. This implies that 
$z\in Z(\cM(x))\mapsto \{w_n\}$ is an injective map. 

Suppose that a bounded harmonic sequence $\{w_n\}$ is given, and $M=\sup_n\|w_n\|_\infty$. 
We regard $w_n\in Z(A_n(x))$. 
Let $C_n:=\{a\in \cM(x)\mid F_n(a)=w_n, \|a\|\leq M\}$. Then $\{C_n\}$ is 
a decreasing sequence of 
non-empty $\sigma$-weakly compact sets, since $w_m\in C_n$ for $m\geq n$. 
Hence $\bigcap_n C_n\ne \emptyset$. Take $a\in \bigcap_n C_n\ne \emptyset$. 
By the martingale convergence, $a=\lim_n w_n$ in $\sigma$-weak topology. Let $u\in A_m(x)$ be an arbitrary unitary.
Then $uau^*=\lim_nuw_nu^*$, and $uw_nu^*=w_n$ for $n\geq m$. Thus $uau^*=a$, and hence 
$a\in Z(\cM(x))$ holds.
These results  implies 
$z\in Z(\cM(x))\mapsto \{w_n\}$ is a bijection map.

Take $\mu$ as in Theorem \ref{thm:groupoidReiter2}. 
Then $w_n^x(\cdot )$ is a constant function on $\cG^x$ by Theorem \ref{thm:groupoidsTailTrivial}, and   
hence $\cM(x)$ is a factor. 
\hfill$\Box$

\bigskip

In the following, we fix a $\mu$ as in Theorem \ref{thm:groupoidReiter2}.
We next consider the lift of $\alpha$ to $\cM(x)$.

\begin{lem}\label{lem:RNderivative}
Define a positive operator $d_g^n=\bigoplus_{k\in \cG^x}d_g^n(k)$ affiliated with $A_n(x)$ by 
\[
 d_g^{n}(k)=
\sum_{\boldsymbol{t}\in \cG^{r(k),n}_{s(k)}}\sqrt{\frac{\sum\limits_{y\sim r(g)}\mu(gkt_n^{-1}\cdots
 t_1^{-1}g^{-1}(r(g),y))\mu((y,r(g))gt_1)}{\sum\limits_{z\sim r(k)}\mu(kt_n^{-1}\cdots t_1^{-1}(r(k),z))
\mu((z,r(k))t_1)}}e_{\boldsymbol{t},\boldsymbol{t}}. 
\]
$(1)$  We have $\psi^{r(g)}_n(\alpha_g(a))=\psi^{s(g)}_n(d_g^nad_g^n)$ for $a\in A_n(s(g))$. \\
$(2)$ $\phi_n^x(d_g^n)=d_g^{n+1}$.
\end{lem}
\textbf{Proof.}
(1) The left hand side is as follows;
\begin{align*}
 \psi^{r(g)}(\alpha_g(a))&=\sum_{k\in \cG^{s(g)}}
\Tr(\rho_n(gk)\alpha_g(a)(k)) \\
&=\sum_{k\in \cG^{s(g)}}\Tr(u_g(k)^*\rho_n(gk)u_g(gk)a(k))\\ &=
\sum_{k\in \cG^{s(g)}} \psi_n^{s(g)}(u_g(k)^*\rho_n(gk)u_g(gk)\rho_n(k)^{-1}a(k)). 
\end{align*}

Here
\begin{align*}
\lefteqn{u_g(gk)^*\rho_n(gk)u_g(gk)\rho_n(k)^{-1}} \\
&=
\sum_{\boldsymbol{t}\in \cG^{r(k),(n)}_{s(k)}}
\sum_{y\sim r(k)}
 {\mu}^{r(k)}(gkt_n^{-1}t_{n-1}^{-1}\cdots t_{1}^{-1}g^{-1}(r(k),y))\mu^{y}((y,r(k))gt_1)\mu^{x_2}(t_2)\cdots \mu^{x_n}(t_n)
e_{\boldsymbol{t},\boldsymbol{t}} \\
&\times 
\left(
\sum_{z\sim r(k)}
{\mu}^{r(k)}(kt_n^{-1}t_{n-1}^{-1}\cdots t_{1}^{-1}(r(k),z))\mu^{z}((z,r(k))t_1)\mu^{x_2}(t_2)\cdots \mu^{x_n}(t_n)
e_{\boldsymbol{t},\boldsymbol{t}}\right)^{-1} \\
 &=d^n_g(k)^2.
\end{align*}
Thus $\psi^{r(g)}_n(\alpha_g(a))=\psi^{s(g)}_n(d_g^nad_g^n)$ holds. (Note that $\rho_n(k)$ and $d_g^n(k)$ commute.) \\
(2)
\begin{align*}
\lefteqn{ \phi_n(d^n_g)(k)} \\
&=\sum_{l\in \cG_{s(k)}}I_{n+1}(kl^{-1},l)(d_g^n(kl^{-1})\otimes 1)I_{n+1}(kl^{-1},l)^* \\
&=\sum_{l\in \cG_{s(k)}}
\left(\sum_{\boldsymbol{t}\in \cG^{r(k),(n)}_{r(l)}}\sqrt{\frac{\sum\limits_{y\in [r(g)]}\mu(gkl^{-1}t_n^{-1}\cdots
 t_1^{-1}g^{-1}(r(g),y))\mu((y,r(g))gt_1)}{\sum\limits_{z\in [r(k)]}\mu(kl^{-1}t_n^{-1}\cdots t_1^{-1}(r(k),z))
\mu((z,r(k))t_1)}}e_{\boldsymbol{t},\boldsymbol{t}}\otimes e_{ll}\right) \\
&=\sum_{\boldsymbol{t}\in \cG^{r(k),(n+1)}_{s(k)}}\sqrt{\frac{\sum\limits_{y\in [r(g)]}\mu(gkt_{n+1}^{-1}t_n^{-1}\cdots
 t_1^{-1}g^{-1}(r(g),y))\mu((y,r(g))gt_1)}{\sum\limits_{z\in [r(k)]}\mu(kt_{n+1}^{-1}t_n^{-1}\cdots t_1^{-1}(r(k),z))
\mu((z,r(k))t_1)}}e_{\boldsymbol{t},\boldsymbol{t}} \\
 &=d_g^{n+1}(k)
\end{align*}

\hfill$\Box$

By Lemma \ref{lem:RNderivative}, there exists a positive operator $d_g^x$ affiliated  with $\cM(x)$
such that $\psi^{r(g)}(\alpha_g (a))=\psi^{s(g)}(d_g^{s(g)}ad_g^{s(g)})$. 

\begin{lem}\label{lem:exttoGNS}
We can extend $\alpha_g$ to an isomorphism in  $\mathrm{Iso}\bigl(\cM(s(g)),\cM(r(g))\bigr)$.
\end{lem}
\textbf{Proof.} 
Let $x=s(g)$, $y=r(g)$.
Define an operator $U\in B(H_x,H_y)$ by
$U(a\xi_x)=\alpha_g(a(d_g^{x})^{-1})\xi_y$. 
Then 
\[
 \|U(a\xi_\psi)\|^2=\|\alpha_g(a(d_g^x)^{-1})\xi_\psi\|^2=\psi^y\left(\alpha_g((d_g^x)^{-1}a^*a(d_g^x)^{-1})\right)
=\psi^x(a^*a)=\|a\xi_x\|^2.
\]
An adjoint 
$U^*$ is given by $U^*(a\xi_y)=\alpha_g^{-1}(a)d_g\xi_x$. Then
it is easy to see $U^*U=UU^*=1$. Hence $U$ is a unitary.
We can verify $U\pi_x(a)U^*=\pi_y(\alpha_g(a))$ as follows:
\begin{align*}
 U\pi_x(a)U^*b\xi_y&=
U\pi_x(a)\alpha_g^{-1}(b)d_g^x\xi_x=
Ua\alpha_g^{-1}(b)d_g^x\xi_x \\ 
&=\alpha_g(a\alpha_g^{-1}(b)d_g^x(d_g^x)^{-1})\xi_y=\pi_y(\alpha_g(a))b\xi_y.
\end{align*}
Similarly, we have $ U^*\pi_y(a)U=\pi_x(\alpha_g^{-1}(a))$. 
Thus $U\cM(x) U^*=\cM(y)$ holds, and $\Ad U$ is an extension of $\alpha_g$. \hfill$\Box$

\begin{thm}\label{thm:modelfreeaction}
Let $\cG$ be an amenable measurable discrete groupoid. 
Then there exists an outer free action  of $\cG$ on the injective factor
 of type II$_1$
 $\cR_0$ with given 3-cocycle $c$.
\end{thm}
\textbf{Proof.} Let $(\alpha,c)$ be an outer action on
$\{\cM(x)\}_{x\in \cG^{(0)}}$
constructed in Lemma \ref{lem:exttoGNS}.
Let $\cR_\infty$ be the injective factor of type III$_1$, and 
consider $\alpha\otimes \id$ on $\cM(x)\otimes \cR_{\infty}$, which is an injective factor of type III$_1$. 
Then the canonical extension $\widetilde{\alpha\otimes\id}$ is a trace preserving outer action on an injective factor of type
II$_\infty$. By the standard argument, we get an outer action
$(\alpha^0,w^0,c)$ of $\cG$
with 3-cocycle $c$ on an injective factor of type II$_1$ $\cR_0$. 

If we apply the above construction to a trivial 3-cocycle, then we get a genuine faithful action of $\cG$ on 
$\cR_0$.
By taking an infinite tensor product, we get a free genuine action $\sigma$ of $\cG$ on $\cR_0$.
(Such an action is also constructed in \cite{Su-Tak-RIMS}.)
Then 
$(\alpha^0\otimes \sigma, w^0\otimes 1, c)$ is a desired one.
\hfill$\Box$

\subsection{Realization of actions with given invariants}\label{subsec:actionInv}
Once we have a free outer action of $\cG$ on $\cR_0$, we can construct
model actions with given invariants as in \cite{M-unif-Crelle}.
Let $(\alpha^0,w^0,c)$ be  a free outer action of $\cG$ with 3-cocycle
$c$ on $\cR_0$.

Let us  $\{\cN, \beta, \chi\}$ 
be given, where $\cN\subset \cG$ is a normal subgroupoid, 
$\beta$ is a  genuine action of $\cG$ on a family of ergodic flows $\{\cC(x),\theta^x\}$
 such that 
$\beta_n=\id$ for $n\in \cN$,  
and 
$\chi=(\lambda,\mu,d)$ is a characteristic invariant.
Thus we are given a family of ergodic flows
$\{\cC(x),\theta^x\}$, and $\chi=(\lambda,\mu,d)$ for some
$\lambda(n,g),\mu(m,n),d(n,t)\in \cU(\cC(x))$.
 We will construct an outer action $\alpha$ with
$\mathrm{Inv}(\alpha)=\{\cN, \beta, \chi\}$. 
 Let  $\cP(x)$ be an injective factor, whose flow
of weights is $\{\cC(x),\theta^x\}$.

At first we recall the following fact \cite{Su-Tak-fields}.
\begin{thm}
 Let $\cR$ be an injective factor, and $\{\cC,\theta\}$ the flow of
 weights of $\cR$. Then the following exact sequence
 splits.
\[
1\longrightarrow \overline{\Int}(\cR)\longrightarrow
 \Aut(\cR)\overset{\md}{\longrightarrow} \Aut_\theta(\cC)\longrightarrow 1.
 \]
\end{thm}

In a similar way, we can show an action $\beta$ of $\cG$ on
$\{\cC(x),\theta^x\}_{x\in \cG^{(0)}}$
can be lift to that on $\{\cP(x)\}_{x\in \cG^{(0)}}$, which we denote by $\beta$.

Let $\varphi_x$ be a dominant weight on $\cP(x)$,
$\cQ(x):=\cP(x)_{\varphi_x}$, $\cP(x)=\cQ(x)\rtimes_{\theta^x}
\mathbb{R}$ the continuous decomposition, and $u^x(s)\in \cP(x)$ be the 
implementing unitary. We may assume that $\varphi_{s(g)}=\varphi_{r(g)}\circ
\beta_g$ for all $g\in \cG$.

Let $\sigma_n^x:=\sigma^{\varphi_x}_{d(n)^*}$ be an extended modular
automorphism for 1-cocycle $d(n,\cdot)^*$. Then
$(\sigma_n^x,\mu(\cdot,\cdot))$ is an outer action of $\cN_x$ with 
3-cocycle $\overline{c(\cdot,\cdot,\cdot)}$ by the relation (CC2) and (CC3). 
Hence $(\sigma^x_m\otimes \alpha^{0}_m,\mu(m,n)\otimes
w^{0}(m,n))$ is a free cocycle crossed action of $\cN_x$ on $\cP(x)\otimes \cR_0$.

Let $\cM(x):=(\cP(x)\otimes \cR_0)\rtimes \cN_x$, and define $\alpha_g\in
\mathrm{Iso}(\cM(x),\cM(y))$, $x=s(g)$, $y=r(g)$ by 
\[
 \alpha_g(a)=\beta_g\otimes \alpha_g^{0}(a),\,\,\, a\in
 \cP(x)\otimes \cR_0,
\]
\[
 \alpha_g(v^x(g^{-1}ng))=(\lambda^y(n,g)\otimes d^{0}(n,g))v^y(n),
\]
where $v^x(n)$ is the implementing unitary in $\cM(x)$, and $d^{0}(n,g)=w^{0}(g,g^{-1}ng)w^{0}(n,g)^*$.

\begin{thm}\label{thm:realizationInv}
Let $\alpha$ be an outer action of $\cG$ constructed above. 
Then we have $\mathrm{Inv}(\alpha)=(\beta, \cN, \chi)$.
\end{thm}

Proof is similar to that of \cite{M-unif-Crelle}, so we omit the detail.

\appendix
\section{Relation with Katayama-Takesaki's classification}\label{sec:KT}

Let $(\alpha, \ual, c)$ be an outer action of a discrete group $G$ on a factor $\cM$. 
We will describe outline of correspondence between $\mathrm{Inv}(\alpha)=(N,\md(\alpha), \chi)$ 
and Katayama-Takesaki's invariant, and explain their classification theorem follows from Theorem \ref{thm:classgroupouter}.  

Let $Q:=G/N$, and $\pi: G\rightarrow  Q$ the quotient map.
In this section, we use letters $g,h,k$ for general elements in $G$, 
$p,q,r$ for those in $Q$, and $l,m,n$ for those in $N$.

Fix a section $p\in Q \rightarrow \secp \in G$, and  
let $\mathfrak{n}(p,q)=\secp\sq\secpq^{-1}\in N$.

We recall the definition of modular obstruction $\mathrm{Ob}_m$ defined by Katayama-Takesaki \cite{KtT-outerI}.
At first note that $\tal_{\secp}$ is a $Q$-kernel on $\tM$. 
Fix $z\in C^2(Q, \cU(\tM))$ such that
$\tal_{\secp}\circ \tal_{\sq}=\Ad z(p,q)\circ \tal_{\secpq}$.
Then we get $d_1\in Z^3(Q, \cU(\cC))$ and $d_2(s;p,q)\in \cU(\cC)$ by  
\[
 \tal_{\secp}(z(q,r))z(p,qr)=d_1(p,q,r)z(p,q)z(pq,r),\,\, \theta_s(z(p,q))=d_2(p,q; s)z(p,q)
\]
We can see 
\[
d_1d_2((p,s), (q,t), (r,u)):=d_1(p,q,r)d_2(s;q,r)\in Z^3(Q\times \mathbb{R},\cC(\cU)) 
\]
and hence get 
an element in $H^3(Q\times \mathbb{R},\cC(\cU))$.
We have another invariant 
\[
 \nu:N\rightarrow H^1_\theta(\mathbb{R},\cU(\cC))
\]
by $\nu(n)=[d(n,\cdot)]$. 
The invariant $\nu$ relates with $d_2$ by 
$\nu(\mathfrak{n}(p,q))=[d_2(\cdot;\mathfrak{n}(p,q))]$ in $H^1_\theta(\mathbb{R},\cU(\cC))$. 
By definition,  $\mathrm{Ob}_{\mathrm{m}}(\alpha)=([d_1d_2],\nu)$ is  the modular obstruction for $\alpha$.

In what follows, we denote the canonical extension of $\alpha_{\secp}$
by the same symbol $\alpha_{\secp}$.

To describe 
 relation between $\mathrm{Inv}(\alpha)=(N,\md(\alpha), \chi)$ and $\mathrm{Ob}_{\mathrm{m}}(\alpha)$, 
it is convenient to replace $(\alpha,\ual)$ with suitable form.  

Note that any element of $G$ can be expressed as $l\secp$, $l\in N$, $p\in Q$,  uniquely .
Define $\hat{\alpha}_{l\secp}=\alpha_l\circ \alpha_{\secp}$. 
Let 
\begin{align*}
\widehat{\ual}(m,n)&=\ual(m,n), \\
\widehat{\ual}(\secp, n)&=\ual(\secp, n)\ual(\secp n \secp^{-1}, \secp)^*, \\
 \widehat{\ual}(\secp, \sq)&=\ual(\secp, \sq)\ual(\mathfrak{n}(p,q),\secpq)^*
\end{align*}
for $m,n\in N$, $p,q\in Q$.
Then we have
\begin{align*}
 \hat{\alpha}_m\circ\hat{\alpha}_n&=\Ad \widehat{\ual}(m,n)\circ \hat{\alpha}_{mn},  \\
 \hat{\alpha}_{\secp}\circ\hat{\alpha}_n&=\Ad \widehat{\ual}(\secp,n)\circ {\alpha}_{\secp n\secp^{-1}}\circ \alpha_{\secp}
=\Ad \widehat{\ual}(\secp,n)\circ \hat{\alpha}_{\secp n},  \\
 \hat{\alpha}_{\secp}\circ\hat{\alpha}_{\sq}&
=\Ad \widehat{\ual}(\secp,\sq)\circ {\alpha}_{\mathfrak{n(p,q)}}\circ \alpha_{\secpq}
=\Ad \widehat{\ual}(\secp,\sq)\circ \hat{\alpha}_{\secp\sq}.
\end{align*}
Define 
\[
 \widehat{\ual}(m\secp, n\sq)=
\alpha_m(\widehat{\ual}(\secp, n))
\alpha_{m}\alpha_{\secp n \secp^{-1}}\left(
\widehat{\ual}(\secp,\sq)\right)
\widehat{\ual}(m,\secp n \secp^{-1})
\widehat{\ual}(m\secp n \secp^{-1},\mathfrak{n}(p,q)).
\]
We can see that $\hat{\alpha}_{g}\circ \hat{\alpha}_h=\Ad
\widehat{\ual}(g,h)\circ \hat{\alpha}_{gh}$. Let $\hat{c}\in Z^3(G,\mathbb{T})$ be
a 3-cocycle associated with $(\hat{\alpha},\widehat{\ual})$. Then we have
 $[\hat{c}]=[c]$ in $H^3(G,\mathbb{T})$. 

In what follows, we replace $(\alpha,\ual)$ with
$(\hat{\alpha},\widehat{\ual})$. In particular, we have $\ual(n,\secp)=1$ and hence
$\alpha_{n\secp}=\alpha_n\alpha_{\secp}$.
We also fix $z(p,q)$ as 
$z(p,q)=\ual(\secp,\sq)\tilu_{\mathfrak{n}(p,q)}$. 
Then we have $\alpha_{\secp}\circ \alpha_{\sq}=\Ad z(p,q)\circ \alpha_{\secpq}$, and 
$d(\mathfrak{n}(p,q),s)=d_2(\mathfrak{n}(p,q;s))$.

\begin{df}\label{df:HJRmap}
 Let $\chi=[\lambda,\mu,d]$ be a characteristic invariant for $(\alpha,\ual)$.
Define $\delta[\lambda,\mu]\in C^3(Q,\cU(\cC))$ by
\[
 \delta[\lambda,\mu](p,q,r)
=\lambda(\secp \mathfrak{n}(q,r)\secp^{-1},\secp)
\mu(\secp\mathfrak{n}(q,r)\secp^{-1},\mathfrak{n}(p,qr))
\mu(\mathfrak{n}(p,q),\mathfrak{n}(pq,r))^*.
\]
\end{df}

The map $\delta$ is an analogue of Huebschmann-Jones-Ratcliffe map, \cite{Hueb}, \cite{J-act}, \cite{Ra}.

At first we express $d_1$ in terms of $\lambda,\mu$ and $c$.
\begin{lem}\label{lem:rel-I}
We have
 \[
  d_1(p,q,r)=c(\secp,\sq,\secr)\delta[\lambda,\mu](p,q,r).
 \]
\end{lem}
\textbf{Proof.} 
By 
\[
 \alpha_{\secp}(\ual(\secp,\sq))\ual(\secp,\sq\secr)=c(\secp, \sq,\secr)\ual(\secp,\sq)\ual(\secp\sq,\secr),
\]
we have 
\begin{align*}
& \alpha_{\secp}(\ual(\secp,\sq))\ual(\secp,\mathfrak{n}(q,r))\alpha_{\secp \mathfrak{n}(q,r)\secp^{-1}}(\ual(\secp,\sqr))
\ual(\secp \mathfrak{n}(p,q)\secp^{-1},\mathfrak{n}(p,qr))\\
&=c(\secp, \sq,\secr)\ual(\secp,\sq)\alpha_{\mathfrak{n}(p,q)}(\ual(\secpq,\secr))
\ual(\mathfrak{n}(p,q),\mathfrak{n}(pq,r)).
\end{align*}

Recall that we have fixed  $z\in C^2(Q,\cU(\cC))$ as
$z(p,q)=\ual(\secp,\sq)\tilu_{\mathfrak{n}(p,q)}
$.
Then
\begin{align*}
 d_1(p,q,r)&=\alpha_{\secp}(z(q,r))z(p, qr)z(pq, r)^*z(p,q)^* \\
&=\alpha_{\secp}\left(\ual(\sq,\secr)\tilu_{\mathfrak{n}(q,r)}\right)
\ual(\secp,\sqr)\tilu_{\mathfrak{n}(p,qr)}
\tilu_{\mathfrak{n}(pq,r)}^*\ual(\secpq,\secr)^*
\tilu_{\mathfrak{n}(p,q)}^*\ual(\secp,\sq)^* \\
&=\lambda(\secp \mathfrak{n}(q,r)\secp^{-1},\secp)
\alpha_{\secp}\left(\ual(\sq,\secr)\right)
\ual(\secp, \mathfrak{n}(q,r))
\tilu_{\secp\mathfrak{n}(q,r)\secp^{-1}}
\ual(\secp,\sqr)\tilu_{\mathfrak{n}(p,qr)} \\
&\times \tilu_{\mathfrak{n}(pq,r)}^*\ual(\secpq,\secr)^*
\tilu_{\mathfrak{n}(p,q)}^*\ual(\secp,\sq)^* \\
&=\lambda(\secp \mathfrak{n}(q,r)\secp^{-1},\secp)
\alpha_{\secp}\left(\ual(\sq,\secr)\right)
\ual(\secp, \mathfrak{n}(q,r))
\alpha_{\secp\mathfrak{n}(q,r)\secp^{-1}}\left(
\ual(\secp,\sqr)
\right) \\
&\times 
\tilu_{\secp\mathfrak{n}(q,r)\secp^{-1}}\tilu_{\mathfrak{n}(p,qr)}
\tilu_{\mathfrak{n}(pq,r)}^*\tilu_{\mathfrak{n}(p,q)}^*
\alpha_{\mathfrak{n}(p,q)}\left(
\ual(\secpq,\secr)^*\right)
\ual(\secp,\sq)^* \\
&=\lambda(\secp \mathfrak{n}(q,r)\secp^{-1},\secp)
\mu(\secp\mathfrak{n}(q,r)\secp^{-1},\mathfrak{n}(p,qr))
\mu(\mathfrak{n}(p,q),\mathfrak{n}(pq,r))^* \\
&\times \alpha_{\secp}\left(\ual(\sq,\secr)\right)
\ual(\secp, \mathfrak{n}(q,r))
\alpha_{\secp\mathfrak{n}(q,r)\secp^{-1}}\left(\ual(\secp,\sqr)\right)
\ual(\secp\mathfrak{n}(q,r)\secp^{-1},\mathfrak{n}(p,qr)) \\
&\times 
\ual(\mathfrak{n}(p,q),\mathfrak{n}(pq,r))^*
\alpha_{\mathfrak{n}(p,q)}\left(
\ual(\secpq,\secr)^*\right)
\ual(\secp,\sq)^* \\
&=c(\secp,\sq,\secr)\delta[\lambda,\mu](p,q,r).
\end{align*}
\hfill$\Box$

\medskip

Conversely, we can recover the 3-cocycle $c$ from $d_1,\lambda,\mu$ as
following lemma.
\begin{lem}\label{lem:rel-II}
Put 
\[
 a(m\secp,n\sq)=
\lambda(\secp n \secp^{-1},\secp)^*
\mu(m,\secp n \secp^{-1})^*
\mu(m\secp n \secp^{-1},\mathfrak{n}(p,q))^*.
\]
Then
\[
 c(g,h,k)=\partial(a)(g,h,k)d_1(\pi(g),\pi(h),\pi(k)).
\] 
Here 
\[
 \partial(a)(g,h,k)=\alpha_g(a(h,k))a(g,hk)a(gh,k)^*a(g,h)^*.
\]
\end{lem}
\textbf{Proof.}
At first we express $\ual(g,h)$ by $z,\lambda,\mu,\tilu_n$.
%
Let $g=m\secp$, $h=n\sq$.
\begin{align*}
 \lefteqn{\ual(m\secp,n\sq)} \\
&=\alpha_m({\ual}(\secp, n))
\alpha_{m}\alpha_{\secp n \secp^{-1}}\left(
\ual(\secp,\sq)\right)
{\ual}(m,\secp n \secp^{-1})
{\ual}(m\secp n \secp^{-1},\mathfrak{n}(p,q)) \\
&=\alpha_m({\ual}(\secp, n))
{\ual}(m,\secp n \secp^{-1})
\alpha_{m\secp n \secp^{-1}}\left(
{\ual}(\secp,\sq)\right)
{\ual}(m\secp n \secp^{-1},\mathfrak{n}(p,q))\\
&=\alpha_m({\ual}(\secp, n))
{\ual}(m,\secp n \secp^{-1})
\tilu_{m\secp n \secp^{-1}}
z(p,q)\tilu_{\mathfrak{n}(p,q)}^*
\tilu_{m\secp n \secp^{-1}}^*
{\ual}(m\secp n \secp^{-1},\mathfrak{n}(p,q)) \\
&=
\mu(m,\secp n \secp^{-1})^*
\mu(m\secp n \secp^{-1},\mathfrak{n}(p,q))^*
\alpha_m({\ual}(\secp, n))
\tilu_{m}
\tilu_{\secp n \secp^{-1}}
z(p,q)
\tilu_{m\secp n \secp^{-1}\mathfrak{n}(p,q)} \\
&=
\mu(m,\secp n \secp^{-1})^*
\mu(m\secp n \secp^{-1},\mathfrak{n}(p,q))^*
\tilu_{m}
{\ual}(\secp, n)
\tilu_{\secp n \secp^{-1}}
z(p,q)
\tilu_{m\secp n \secp^{-1}\mathfrak{n}(p,q)} \\
&=
\mu(m,\secp n \secp^{-1})^*
\mu(m\secp n \secp^{-1},\mathfrak{n}(p,q))^*
\lambda(\secp n \secp^{-1},\secp)^*
\tilu_m
\alpha_{\secp}(\tilu_{n})
z(p,q)
\tilu_{m\secp n \secp^{-1}\mathfrak{n}(p,q)} 
\end{align*}

Put 
\[
 a(m\secp,n\sq)=
\lambda(\secp n \secp^{-1},\secp)^*
\mu(m,\secp n \secp^{-1})^*
\mu(m\secp n \secp^{-1},\mathfrak{n}(p,q))^*.
\]

For $g\in G$, let $p=\pi(g)$ and $n(g)=g\secp^{-1}$. Then we have 
\[
 \ual(g,h)=
a(g,h)\tilu_{n(g)}
\alpha_{\secp}(\tilu_{n(h) })
z(\pi(g),\pi(h))\tilu_{n(gh)}^*. 
\]

We compute the
3-cocycle  $c$ for $(\alpha,\ual)$. Let $p=\pi(g)$, $q=\pi(h)$, $r=\pi(k)$.
Then 
\begin{align*}
\lefteqn{c(g,h,k)} \\
 &=
 \alpha_{g}(\ual(h,k))\ual(g,hk)\ual(gh,k)^*\ual(g,h)^* \\
&=
\alpha_{n(g)}\alpha_{\secp}\left(a(h,k)\tilu_{n(h)}\alpha_{\sq}(\tilu_{n(k) })
z(q,r)\tilu_{n(hk)}^*\right)  \times a(g,hk)\tilu_{n(g)}
\alpha_{\secp}(\tilu_{n(hk) })
z(p,qr)\tilu_{n(ghk)}^* \\ 
&\hspace{12pt} \times\left( a(gh,k)\tilu_{n(gh)}
\alpha_{\secpq}(\tilu_{n(k) })z(pq,r)\tilu_{n(ghk)}^*  \right)^* 
\times\left( a(g,h)\tilu_{n(g)}
\alpha_{\secp}(\tilu_{n(h) })z(p,q)\tilu_{n(gh)}^*  \right)^* \\
&=\partial (a)(g,h,k)
\alpha_{n(g)}\alpha_{\secp}\left(\tilu_{n(h)}\alpha_{\sq}(\tilu_{n(k) })
z(q,r)\tilu_{n(hk)}^*\right) 
 \times \tilu_{n(g)}
\alpha_{\secp}(\tilu_{n(hk) })
z(p,qr)\tilu_{n(ghk)}^*  \\
&\hspace{12pt}\times
\tilu_{n(ghk)}z(pq,r)^*\alpha_{\secpq}(\tilu_{n(k) }^*)\tilu_{n(gh)}^*
\times\tilu_{n(gh)}z(p,q)^*\alpha_{\secp}(\tilu_{n(h) }^*)\tilu_{n(g)}^*
\\
&=\partial (a)(g,h,k)
\tilu_{n(g)}\alpha_{\secp}\left(\tilu_{n(h)}\alpha_{\sq}(\tilu_{n(k) })
z(q,r)\tilu_{n(hk)}^*\right) \\
& \hspace{12pt}\times 
\alpha_{\secp}(\tilu_{n(hk) })
z(p,qr) 
z(pq,r)^*\alpha_{\secpq}(\tilu_{n(k) }^*)
z(p,q)^*\alpha_{\secp}(\tilu_{n(h) }^*)\tilu_{n(g)}^*
\\
&=\partial (a)(g,h,k)
\tilu_{n(g)}\alpha_{\secp}\left(\tilu_{n(h)}\alpha_{\sq}(\tilu_{n(k) })
\right) \\
&\hspace{12pt}  \times 
\alpha_{\secp}(z(q,r))z(p,qr) 
z(pq,r)^*z(p,q)^*
\alpha_{\secp}\alpha_{\sq}(\tilu_{n(k) }^*)
\alpha_{\secp}(\tilu_{n(h) }^*)\tilu_{n(g)}^*
\\
&=\partial (a)(g,h,k)d_1(p,q,r).
\end{align*}
Thus $c$ is given by 
$c(g,h,k)=\partial (a)(g,h,k)d_1(p,q,r)$.  \hfill$\Box$

\medskip

\noindent
\textbf{Remark.} \\
(1) The formula  in Lemma \ref{lem:rel-II} is obtained by Katayama-Takesaki in \cite[Lemma 2.11]{KtT-outerI}. \\
(2) In Katayama-Takesaki's formulation, $\lambda$ and $\mu$   do not appear explicitly. However we can recover these by
$\nu(mn)=\nu(m)\nu(n)$ and $\md(\alpha_g)(\nu(g^{-1}ng))=\nu(n)$, which are equivalent to  the relation (CC1) and (CC2),
respectively.

\medskip

With a bit of effort, we can show that the above correspondence
preserves equivalence classes of invariants. 
By Theorem \ref{thm:classgroupouter}, we can show Katayama-Takesaki's classification. Namely, we can conclude the following. 
\begin{thm}
 Let $\alpha$ and $\beta$ be outer actions on an injective factor $\cM$. Then $\alpha$ and $\beta$ are conjugate in 
$\mathrm{Out}(\cM)$ if and only if $\sigma(\md(\alpha),
 \mathrm{Ob}_{\mathrm{m}}(\alpha))=(\md(\beta),
\mathrm{Ob}_{\mathrm{m}}(\beta))$
for some $\sigma\in \Aut_\theta(\cC)$.   
\end{thm}


\ifx\undefined\bysame
\newcommand{\bysame}{\leavevmode\hbox to3em{\hrulefill}\,}
\fi

\end{document}